\documentclass[12pt]{article}
\usepackage{amssymb,amsmath, float, txfonts}%, t1enc, txfonts, multirow}%, graphicx,float,epsf,t1enc, txfonts, multirow}
\begin{document}

\def\abstractname{\bf Abstract}
\def\dfrac{\displaystyle\frac}
\let\oldsection\section
\renewcommand\section{\setcounter{equation}{0}\oldsection}
\renewcommand\thesection{\arabic{section}}
\renewcommand\theequation{\thesection.\arabic{equation}}
\newtheorem{theorem}{\indent Theorem}[section]
\newtheorem{lemma}{\indent Lemma}[section]
\newtheorem{proposition}{\indent Proposition}[section]
\newtheorem{definition}{\indent Definition}[section]
\newtheorem{remark}{\indent Remark}[section]
\newtheorem{corollary}{\indent Corollary}[section]
\def\pd#1#2{\displaystyle\frac{\partial#1}{\partial#2}}
\def\d#1{\displaystyle\frac{d#1}{dt}}

\title{\LARGE\bf
Classification of self-similar singular solutions with large mass for Keller-Segel model with signal consumption% and singular sensitivity
%Self-similar singular solutions and instantaneous smoothness  of solutions for Keller-Segel model
%Classification of self-similar singular solutions of Keller-Segel model
\thanks{This work is supported by NSFC(12271186,12171166)}
\\
\author{Chunhua Jin\thanks{
Corresponding author. Email: {\tt jinchhua@126.com}}
\\
\small (\it{School of Mathematical Sciences, South China
Normal University, }
\\
\small \it{Guangzhou, 510631, China}\rm)
}}

\date{}

\maketitle

%Mathematics Subject Classification (2020): 35C06 35A21 35A08

%It is widely recognized that scaling invariance is  a fundamental requirement for the existence of self-similar structural solutions, thereby necessitating the fulfillment of the condition $\kappa(1-\alpha) =1-\frac N2$.
%\begin{align*}
%\left\{
%\begin{aligned}
%&u_t=D_u\Delta u-\chi\nabla\cdot(\frac{u}{v}\nabla v),
%\\
%&v_t=D_v\Delta v-v^{\alpha}u.
%\end{aligned}\right.
%\end{align*}
%It is well-acknowledged that scaling invariance is  a fundamental requirement for the existence of self-similar structural solutions, this necessitates the fulfillment of the condition  $\kappa(1-\alpha) =1-\frac N2$.
%Building on this,
\begin{abstract}
In this paper,  we concentrate on investigating the self-similar singular solutions of Keller-Segel model with signal consumption ($-uv^{\alpha}$) and singular sensitivity.
We perform a detailed exploration into the existence and decay rate of self-similar solutions,
particularly, the permissibility of arbitrary mass for these solutions across all possible cases.
Based on these findings, we can delve deeper into verifying that these self-similar solutions $(u, v)$  exhibit varying degrees of singularity depending on the value of $\alpha$ and the spatial dimension.

Our analysis reveals that the component  $u$ (with arbitrary mass) of the solution consistently behaves analogous to heat kernel, that is, $u$  exhibiting a Dirac $\delta$ initial singularity identical to that of the fundamental solution, and converges to $0$ in the sense of the $L^p$-norm ($p>1$) as time approaches infinity.
However, the initial behavior of the other component  $v$ varies significantly based on the value of $\alpha$ and the spatial dimension, exhibiting regularity  (not singular), less singularity, or strong singularity (more singular than fundamental solution).    Moreover,  both $u$ and $v$ undergo instantaneous smoothing, becoming smooth immediately after $t>0$. This phenomenon  reveals the adaptive strategies of cells in high-density aggregation environments to prevent resource depletion, reflecting an optimization process of self-organizing behavior.
\end{abstract}


\section{Introduction}

The phenomenon of a chemotactic system blowing-up at a finite time can be used to characterize the chemotactic collapse phenomenon in biology, such as the aggregation of slime mold cells to form spores \cite{Na}. The reason for this phenomenon is due to the effects of chemotaxis, organisms
rapidly gather in a specific area, leading to a sharp increase in biological density within that region.
This phenomenon is manifested in mathematical models as rapid changes in certain variables, which
may result in the loss of regularity in the solution. The chemotactic collapse phenomenon is an important regulatory mechanism in the process of slime mold spore formation, ensuring that spores can form at the correct time and location, and be released under appropriate conditions to ensure the reproduction and survival of slime molds. In mathematics, chemotaxis collapse specifically refers to the formation of Dirac $\delta$-singularity at a finite time.

The study of the singularity of solutions is of great significance in both mathematics and biology. In the past few decades, extensive research has been conducted on the finite time blowing-up phenomenon of Keller-Segel models with signal production, and rich results have been achieved.
But what happens to the solution after it blows-up? Will it maintain its singularity, become regularized, or manifest some other phenomenon? This interesting question has long been a focus of researchers. For instance, in 2011, Vazquez et al. \cite{VW} studied the time evolution of an initial singularity  for porous medium equation with fast diffusion and discovered that the initial singularity
can lead to a diverse range of phenomena, such as the solution  maintaining its singularity, or instantaneous smoothing, or regularization after a finite time or infinite time etc. In fact, the study of singular initial value problems for linear, or nonlinear evolution equations has always been a hot topic of concern for mathematical researchers. As early as 1952, Barenblatt \cite{B} proposed the problem of singular solutions. Afterwards, many mathematicians, such as Friedman, Kamin, Peletier, Vazquez, Zuazua et.al \cite{BF, EVZ, K, KP} have successively studied the singular initial value problem of reaction diffusion equations.

Due to their crucial role in the research on the singularity of solutions, the study of self-similar solutions have attracted increasing attention, especially concerning spatial and instantaneous structures, localization, and stability issues. Since Leray \cite{LJ} proposed the problem of self-similar solutions for the 3-D  Navier-Stokes equations, research in this area has flourished. The investigation of self-similar solutions is profoundly significant both physically and geometrically. The structure and existence of these solutions reflect the intricate dynamics resulting from the interaction of various nonlinear and diffusion terms on the properties of the solution.
Such solutions also play an important role in the research into the regularity of nonlinear problems, this has been manifested in the
regularity theory of  harmonic maps and minimal surfaces \cite{TX}.
Self-similar solutions may exhibit precise singular behavior, either at the blowing up time, or at starting time, so the research on self-similar solutions yields clues for studying the singularity of solutions,
and numerous classical results have been achieved in the study of self-similar singular solutions for heat equations with sources or sinks, as well as various nonlinear diffusion equations \cite{BPT, L, PW}. For instance, Brezis and Friedman \cite{BF} demonstrated in 1983 the existence of singular solutions satisfying the condition $u(x, 0)=M\delta(x)$ for the semilinear heat equation
$$
u_t-\Delta u+u^p=0
$$
when $p\in(1, 1+\frac 2N)$.  In the context of linear theory, such a solution is commonly referred to as a fundamental solution, and occasionally, it is also labeled as a "source-type solution" \cite{EVZ}. Subsequently, Brezis, Peletier, and Terman \cite{BPT} discovered a singular solution that exhibited even more singularity than the fundamental solution, which is called "very singular solution".
Furthermore, numerous studies have observed the existence of "very singular solutions" in  numerous equations involving $p$-Laplace diffusion or porous medium diffusion. For more details, please refer to \cite{CQW, KP, KPV, KV, PT, PW}, etc.

Based on numerous results regarding the existence of self-similar solutions mentioned above, the asymptotic behavior of solutions at large time scales has also attracted widespread attention from researchers. The findings indicate that the convergence of solutions is influenced not only by the diffusion term and the reaction term but also by the initial values, particularly the decay rate of the initial values at infinity. These factors will lead to rich and interesting long-time dynamic behaviors. For instance, under different conditions, the solutions may converge to similar solutions of the heat equation , or various similar solutions of the full equation. For details, please refer to the literature \cite{GA, KP1, NY, ZJ}.

The study of self-similar solutions to the Keller-Segel model has been ongoing for over two decades. At the turn of the last century,  Biler, Mizutani, and  Nagai et al. \cite{BCD, BP, BKW, MM, MN} investigated the existence of self-similar solutions in the form of $t^{-1}G_\alpha\left(\frac{|x|}{\sqrt t}\right)$ for the following 2D Keller-Segel system
\begin{align}
\label{1-1}\left\{
\begin{aligned}
&u_t=\Delta u-\chi\nabla\cdot(u\nabla v)
\\
&\tau v_t=\Delta v+u
\end{aligned}\right.
\end{align}
with $\tau=0$, or $\tau>0$.  Notably, when $\tau=0$, it has been established that these self-similar solutions serve as global attractors \cite{BDP, BM1}; when $\tau>0$,  global self-similar solution with arbitrarily large mass has been successfully constructed  for large $\tau$ by Biler,  Corrias, et al. in two dimensional space \cite{BCD}, this also indicates the significant result  that the Keller-Segel model with  parabolic-parabolic  structure  does not necessarily  lead to blow up even when the initial mass exceeds the critical mass.
Recently, Biler et al.  \cite{BKW} also investigated the self-similar solutions in several space dimensions.
Additionally, the backward self-similar blow-up solutions of the Keller-Segel model have garnered attention and been explored by some authors, as detailed in \cite{HMV, LB}.

In addition to the Keller-Segel model with signal production, there have been few studies on the singularity of solutions for Keller-Segel model with signal consumption mechanisms. In 2021, Carrillo, et al.\cite{CLW} explored the boundary spike-layer solution of the following model
\begin{align*}
\left\{
\begin{aligned}
&u_t=u_{xx}-\chi \left(\frac{u}{v}v_x\right)_x
\\
&v_t=\varepsilon v_{xx}-v^\alpha u
\end{aligned}\right.
\end{align*}
in a one-dimensional half-space interval $(0,+\infty)$. They obtained a non-constant steady state $(U, W)$, where $U$ develops into
a Dirac $\delta$-singularity at the boundary $x=0$ as $\chi\to\infty$ or $\varepsilon\to 0$, while $W$ forms a boundary-layer profile as $\varepsilon\to 0$.
 The potential occurrence of finite-time blow-up phenomena in this model has also attracted significant attention. In 1995, Rascle and Ziti \cite{RZ} investigated the self-similar blow-up solutions of the following model
\begin{align*}
\left\{
\begin{aligned}
&u_t=\mu \Delta u-\chi\nabla\cdot(\frac{u}{v}\nabla v)
\\
&v_t=-\kappa v^\alpha u,
\end{aligned}\right.
\end{align*}
and provided  analytical expressions of self-similar blow-up solutions. Although they neglected the diffusion effect of oxygen, it provided valuable insights into the characteristics of blowup solutions. Recently, we have considered the backward self-similar solutions of the model \eqref{1-1} with  non-zero diffusion coefficients \cite{J1}. For the one-dimensional case, we have verified the occurrence of chemotactic collapse, where the solution forms a Dirac $\delta$-singularity at a finite time. However, in the higher-dimensional case, the backward self-similar solutions do not possess finite mass, exhibiting blow-up behavior everywhere.

However, after the solution blows up, how can we describe its subsequent behavior? Will it maintain its singularity, become regularized, or exhibit some other complex phenomena? To explore this question, this paper is dedicated to studying the forward self-similar singular solutions with initial singularities. We will delve into the properties of such singular solutions and the possible evolutionary paths they may follow afterward. That is, we consider the following Keller-Segel model with singular sensitivity and signal consumption
\begin{align}
\label{1-1}\left\{
\begin{aligned}
&u_t=D_u\Delta u-\chi\nabla\cdot\left(\frac{u}{v}\nabla v\right),&& (x,t)\in  \mathbb R^N\times \mathbb R^+,
\\
&v_t=D_v\Delta v-v^{\alpha}u, &&(x,t)\in \mathbb R^N\times \mathbb R^+,
\end{aligned}\right.
\end{align}
where $N\ge 1$, $u$ and $v$ represent the density of bacteria and the concentration of oxygen respectively.
The diffusion coefficients $D_u$ and $D_v$  are positive constants, $\chi\nabla\cdot(\frac{u}{v}\nabla v)$
represents the chemotaxis term, the positive constant $\chi$  denotes the chemotaxis sensitivity coefficient, $v^{\alpha}u$ signifies the consumption of oxygen by bacteria with $\alpha\ge 0$.
This model was proposed in 1971 by Keller and Segel \cite{KS} to  describe the wave phenomenon discovered in Adler's experiment \cite{A1, A2}.

It is well-acknowledged that scaling invariance is  a fundamental requirement for the existence of self-similar structural solutions.
Noticing that if $(u(x,t),v(x,t))$ solves \eqref{1-1}, so do the rescaled functions
$$
(u_\rho(x,t), v_\rho(x,t))=(\rho^Nu(\rho x, \rho^2t),  \rho^{-2\kappa}v(\rho x, \rho^2t)),
$$
when $(1-\alpha)\kappa=1-\frac{N}2$.  Therefore, we assume that
\begin{equation}\label{1-3}
\kappa(1-\alpha)=1-\frac{N}2.
\end{equation}
As a result, we can explore the self-similar solution of \eqref{1-1} by taking $\rho=\frac1{\sqrt t}$. Specifically, let
\begin{equation}\label{1-2}
u(x,t)=t^{-\frac{N}2}\varphi(r), \qquad v(x,t)=t^{\kappa}\psi(r),
\end{equation}
with $r=t^{-\frac12}|x|$.  Substituting these transformations into the system \eqref{1-1} leads to
\begin{align}\label{1-4}
\left\{
\begin{aligned}
&-\frac{N}2\varphi-\frac12r\varphi'(r)=D_u\left(\varphi''(r)+\frac{N-1}{r}\varphi'(r)\right)-\chi\left(\left(\frac{\varphi}{\psi} \psi'(r)\right)'+\frac{N-1}{r}\frac{\varphi}{\psi} \psi'(r)\right),
\\
&\kappa\psi-\frac{1}2r \psi'(r)=D_v\left(\psi''(r)+\frac{N-1}{r}\psi'(r)\right)-\varphi\psi^\alpha.
\end{aligned}\right.
\end{align}
Given the radial symmetry property, it is straightforward to deduce that $\varphi'(0)=\psi'(0)=0$. Subsequently, we present the initial condition as follows:
\begin{align}\label{1-5}
\varphi(0)=A, \quad \psi(0)=B, \quad \varphi'(0)=0, \quad \psi'(0)=0,
\end{align}
where $A$, $B$ are positive constants to be determined.
From the first equation of \eqref{1-4}, we derive that
$$
-\frac{1}2 (r^N\varphi)'=D_u(r^{N-1}\varphi')'-\chi\left(r^{N-1}\frac{\varphi}{\psi} \psi'(r)\right)'.
$$
Integrating the above equation form $0$ to $r$ yields
$$
-\frac{1}2r\varphi=D_u\varphi'-\chi\frac{\varphi}{\psi} \psi'(r),
$$
that is
$$
\left(D_u\ln\varphi-\chi\ln\psi+\frac14r^2\right)'=0.
$$
Hence,  we arrive at
\begin{equation}\label{1-6}
\varphi(r)=AB^{-\frac{\chi}{D_u}}\psi^{\frac{\chi}{D_u}}e^{-\frac{1}{4D_u}r^2}.
\end{equation}
Therefore, the initial value problem of system \eqref{1-4}   can be converted into the initial value problem of the following single equation:
\begin{align}\label{1-7}
\left\{
\begin{aligned}
&\kappa \psi-\frac{1}2r \psi'(r)=D_v\left(\psi''(r)+\frac{N-1}{r}\psi'(r)\right)- AB^{-\frac{\chi}{D_u}} \psi^{\alpha+\frac{\chi }{D_u}}e^{-\frac{1}{4D_u}r^2},
\\
&\psi(0)=B, \quad \psi'(0)=0.
\end{aligned}\right.
\end{align}
By setting  $\phi=\frac{\psi}{B}$, the problem \eqref{1-7} is  ultimately reduced to
\begin{align}\label{1-8}
\left\{
\begin{aligned}
&\kappa \phi-\frac{1}2r\phi'(r)=D_v\left(\phi''(r)+\frac{N-1}{r}\phi'(r)\right)-A B^{\alpha-1} \phi^{\alpha+\frac{\chi}{D_u}}e^{-\frac{1}{4D_u}r^2},
\\
&\phi(0)=1, \quad \phi'(0)=0,
\end{aligned}\right.
\end{align}
where $(1-\alpha)\kappa=1-\frac{N}2$.  It can also be expressed in the equivalent form as follows:
\begin{align}\label{1-9}
\left\{
\begin{aligned}
&D_v\left(\rho\phi'\right)'=\kappa \rho\phi+A B^{\alpha-1}\rho e^{-\frac{r^2}{4D_u}} \phi^{\alpha+\frac{\chi}{D_u}},
\\
&\phi(0)=1, \quad \phi'(0)=0,
\end{aligned}\right.
\end{align}
where
$$
\rho(r)=r^{N-1}e^{\frac{r^2}{4D_v}}.
$$

In what follows,  our objective is to establish the global existence of solution to Problem \eqref{1-8} (or equivalently, to consider Problem \eqref{1-9}) and to calculate its decay rate at infinity.
We observe that there exists an energy functional as follows
$$
\mathcal F(\phi)=\frac{D_v}2\int_0^\infty\rho|\phi'|^2dr+\frac{\kappa}2\int_0^\infty\rho|\phi|^2dr+\lambda
\int_0^\infty\rho\exp\left( -\frac{r^2}{4D_u} \right)|\phi|^{\alpha+\frac{\chi}{D_u}+1}dr
$$
associated with Problem \eqref{1-8}, and the solution corresponds to a critical point of this energy functional.

(i)When  $\kappa< -\frac{N}2$, we discovered that the variational method can be employed to establish the global existence of positive solutions in a weighted Sobolev space.

(ii) However, it seems that the variational method is inapplicable when $\kappa\ge -\frac{N}2$, since $\kappa= -\frac{N}2$ is precisely the principle eigenvalue of the linearized problem from Lemma \ref{lem2-9}. This implies that
$$
\frac{D_v}2\int_0^\infty\rho|\phi'|^2dr+\frac{\kappa}2\int_0^\infty\rho|\phi|^2dr\ge 0, \ \text{for any}\ \phi\in W_{\rho}^{1,2}.
$$
Consequently, there is no positive solution within the aforementioned weighted Sobolev space $W_{\rho}^{1,2}$. In fact,
if $\phi\in W_{\rho}^{1,2}$ is a solution of the problem \eqref{1-8},  it is straightforward to derive that:
$$
0=D_v \int_0^\infty\rho|\phi'|^2dr+\kappa\int_0^\infty\rho|\phi|^2dr+A B^{\alpha-1}\int_0^\infty\rho e^{-\frac{r^2}{4D_u}} \phi^{\alpha+\frac{\chi}{D_u}+1}dx>A B^{\alpha-1}\int_0^\infty\rho e^{-\frac{r^2}{4D_u}} \phi^{\alpha+\frac{\chi}{D_u}+1}dx ,
$$
which indicates that $\phi\equiv 0$. Therefore, the variational method is indeed not applicable to the case $\kappa\ge -\frac{N}2$,
as there is no solution within the weighted Sobolev space $W_{\rho}^{1,2}$.

Nevertheless,  the problem \eqref{1-8}  does admit solutions  that do not belong to this weighted space when $\kappa\ge -\frac{N}2$.
In such cases, we rely on analytical methods, progressively verifying the existence of the solution to Problem \eqref{1-8}.  Additionally, we further investigate  the  decay rate of the global solution at infinity, discovering that depending on the positivity or negativity of $\kappa$, it exhibits either algebraic growth or decay rates, respectively.

 The specific details regarding the global existence of the solution and the estimation of the decay rate are for the case $\kappa\ge -\frac{N}2$ can be stated as follows:

\begin{theorem}\label{thm2-1}
Assume that  $\alpha+\frac{\chi}{D_u}\ge 1$, $\kappa\ge -\frac{N}2$. Then when $\alpha+\frac{\chi}{D_u}=1$, the problem \eqref{1-8} admits a unique global positive solution $\phi\in C^2[0,+\infty)$ for any positive $A, B$; when
 $\alpha+\frac{\chi}{D_u}> 1$, there exists $A^*(\kappa)>0$, such that
when $AB^{\alpha-1}\le A^*$, the problem \eqref{1-8} admits a unique global positive solution $\phi\in C^2[0,+\infty)$,
and when  $AB^{\alpha-1}>A^*$, the problem \eqref{1-8} admits a unique local solution $\phi\in C^2[0, R_{\max})$ with $R_{\max}<\infty$, such that
the solution $\phi(r)$ blows up at $R_{\max}$.
\end{theorem}

\begin{theorem}\label{thm2-2}
Assume that $\alpha+\frac{\chi}{D_u}\ge 1$, and  $\kappa\ge -\frac{N}2$.
Let $\phi\in C^2[0,\infty)$ be the global positive solution of Problem \eqref{1-8}.
Then for both of the following two cases:

i) When $\alpha+\frac{\chi}{D_u}>1$, and $AB^{\alpha-1}$ is appropriately small;

ii) When $\alpha+\frac{\chi}{D_u}=1$, for any $A, B>0$,

\noindent we have
$$
\lim_{r\to\infty}\frac{\phi(r)}{r^{2\kappa}}=M^*,
$$
where $M^*$ is a positive constant.
\end{theorem}

When $\kappa< -\frac{N}2$, it is observed that the solution does not exist globally in the case where $AB^{\alpha-1}$
is either excessively small or excessively large. For this case, the variational approach can be used to derive the global solvability and decay rate estimates for the problem \eqref{1-8}.
\begin{theorem}\label{thm2-3}
Assume that $\kappa< -\frac{N}2$,  $1<\alpha+\frac{\chi}{D_u}<\frac{N+2}{(N-2)_+}$. The solution does not exist globally in the case where $AB^{\alpha-1}$  is either excessively small or excessively large. However, there exists a particular value $A^*>0$, such that
when $AB^{\alpha-1}=A^*$, the problem \eqref{1-8} admits a unique global positive solution $\phi\in C^2[0,+\infty) \cap W_{\rho}^{1,2}(0,\infty)$,
and $\phi(r)=o(\rho^{-\frac12})$ as $r\to\infty$.
Here
$$
W_{\rho}^{1,2}(0,\infty)=\left\{\phi;  \int_0^{\infty}\rho|\phi|^2dr+ \int_0^{\infty}\rho|\phi'|^2dr<\infty\right\}.
$$
\end{theorem}

From the above three Theorems, we can delve deeper into verifying that these self-similar solutions $(u, v)$ of the system \eqref{1-1} exhibit varying degrees of singularity depending on the value of $\kappa$.

Using the aforementioned definitions, if a solution behaves essentially like a fundamental solution at the origin, we refer to it as a "singular solution". If a solution has more singularity than the fundamental solution, it is called a "very singular solution"; On the contrary,
if the singularity of a solution is less  singular than that of the fundamental solution, it is labeled as a "less singular solution".

In fact, our analysis reveals that the element $u$ in the self-similar solution consistently demonstrates Dirac $\delta$ initial singularity  identical to that of the fundamental solution for all $\kappa\in \mathbb R$,   that is
$$
\int_{\mathbb R^N} u(x,t)\phi(x)dx\to M\int_{\mathbb R^N} \phi(x)dx, \quad \text{as } \ t\to 0^+
$$
for any $\phi\in C^{\infty}(\mathbb R^N)\cap L^\infty(\mathbb R^N)$. indicating that $u$ is consistently singular.
However, the behavior of the other element $v$ varies significantly based on the value of $\kappa$. To be precise:

(i) When $\kappa\ge 0$, $v$ is regular, that is $v$ is not singular.

(ii) When $-\frac{N}{2}< \kappa<0$,  $v$ is less singular.

(iii) When $\kappa\le -\frac{N}{2}$, $v$ is very singular.
\\
On the other hand,
both $u$ and $v$ will  become smooth and bounded immediately after $t>0$ (instantaneous smoothing after the chemotaxis collapse occurs). The detailed results described as follows.

 \begin{theorem}\label{thm2-4}
 Assume that $\alpha\not=1$ when $N\not=2$.
 Let $(u, v)$ be the self-similar solution of the model \eqref{1-1} obtained in the above three theorems, which is defined as in \eqref{1-2}. Then  the element $u$  consistently demonstrates Dirac-$\delta$ initial singularity identical to that of the fundamental solution for all values of $\alpha$ and $N$,
$$
u(x,t)\to M\delta(x), \quad \text{as}\ t\to 0^+
$$
in the sense of distribution, where the mass $M$ is defined as
$$
M:=\int_{\mathbb R^N}u(x,t)dx,
$$
which can be arbitrary (see Remark \ref{re1-1}). And for any $1<p\le \infty$,
\begin{equation}\label{1-11}
\lim_{t\to\infty}t^{\frac{N}2(1-\frac1p)}\|u(\cdot, t)\|_{L^p}=M_p,
\end{equation}
where $M_p$ are positive constants depending only on $p$.

However, the initial singularity of $v$ varies with $\alpha$ and $N$,
we give its classification as follows:
\begin{table}[H]
\centering
\begin{tabular}{|c|c|c|c|}
\hline
$N$ & $\alpha$ & Initial behavior of $v$ & Long time behavior of $v$\\
\hline
 & $0 \le \alpha < 1$ & regular & \\
\cline{2-4}
1 & $1 < \alpha \le 2$ & very singular & decay to $0$ at an algebraic rate     \\
\cline{2-4}
  & $\alpha > 2$ & less singular  &    decay to $0$ at an algebraic rate  \\
\hline
 & $\alpha = 1$ & any (regular, less singular, very singular) & depends on initial behavior \\
\cline{2-4}
2 & $\alpha \neq 1$ & regular&   \\
\hline
  & $\alpha > 1$ & regular &  \\
\cline{2-4}
$N\ge 3$ & $0 \le \alpha < \frac{2}{N}$ & less singular &  decay to $0$ at an algebraic rate  \\
\cline{2-4}
  & $\frac{2}{N} \le \alpha < 1$ & very singular &  decay to $0$ at an algebraic rate\\
\hline
\end{tabular}
\caption{ $v$'s behavior based on $N$ and $\alpha$}
\end{table}
\noindent For more details about the singularity of $v$, long time asymptotic behavior, and decay rate, please see  Proposition \ref{pro-2}.

Additionally, $u$ will  become smooth and bounded immediately after $t>0$.
Similarly, for $v$ with initial singularity, it also becomes smooth and bounded immediately after time  $t>0$.
\end{theorem}

Based on the results obtained above, we can see that, after the occurrence of a chemotactic collapse, the solution quickly regularizes, indicating that the system can return to a stable state where variable fluctuations level off. It reveals the adaptive strategies of cells in high-density aggregation environments to prevent resource depletion, where cells avoid over-aggregation by migrating and dispersing, reflecting an optimization process of self-organizing behavior. While, when time is sufficiently large, the self-similar solution $u$ behaves like a heat kernel again, gradually tending towards a spatially uniform distribution state.
In summary, chemotaxis is a fundamental response of organisms to chemical stimuli, influencing their migration and distribution. The chemotactic collapse is an extreme manifestation of this behavior under specific conditions, while the rapid regularization of the solution is a result of the system's self-regulation.

\begin{remark}\label{re1-1}
We observe that the existence of the solution $\phi$ depends on the value of $AB^{\alpha-1}$. Additionally, we remark that the mass $M$ of $u$, given by
$$
M=\int_{\mathbb R^N}u(x,t)dx= A|\partial B_1|\int_0^\infty |\phi(r)|^{\frac{\chi}{D_u}}e^{-\frac{1}{4D_u}r^2}dr
$$
can be arbitrary due to the fact that $A$ can be adjusted arbitrarily by altering the value of $B$.
{\it\bfseries This conclusion indicates that the system \eqref{1-1} admits global solutions with arbitrary mass, and the large mass does not necessarily lead to  the blow-up of solutions at a finite time.}
\end{remark}

\section{Existence of self-similar solutions for the case $\kappa\ge -\frac{N}2$}

As stated in Section 1, throughout the subsequent discussion, we always assume that $A, B>0$.
It is noteworthy that despite the equation \eqref{1-9} admits an energy functional of the following form,
$$
\mathcal F(\phi)=\frac{D_v}2\int_0^\infty\rho|\phi'|^2dr+\frac{\kappa}2\int_0^\infty\rho|\phi|^2dr+\lambda
\int_0^\infty\rho\exp\left( -\frac{r^2}{4D_u} \right)|\phi|^{\alpha+\frac{\chi}{D_u}+1}dr,
$$
the critical point theory cannot be applied to establish the existence of solutions in the space $W_{\rho}^{1,2}$ for the case where $\kappa\ge -\frac{N}2$.  In fact, subsequent analysis reveals that the solutions do not indeed belong to the space $W_{\rho}^{1,2}$. Consequently, we must rely on analytical methods to gradually determine the existence of these solutions.

Employing the standard contraction fixed point method, it becomes straightforward to establish the local existence and uniqueness of the positive solution. Furthermore, utilizing the continuation theorem, we can extend the solution if it remains positive and bounded. Let $[0, R_{\max})$  represent the maximum existence interval for the positive solution. Additionally, from \eqref{1-9}, it becomes evident that when  $\kappa\ge 0$, $\phi$ is increasing on the interval $[0, R_{\max})$. Consequently, we arrive at the following lemma.
\begin{lemma}\label{lem2-1}
Assume that $\kappa\in \mathbb R$, $\alpha+\frac{\chi}{D_u}\ge 1$. Then  the problem \eqref{1-9} admits a unique positive solution $\phi(r)\in C^2[0, R_{\max})$ with $R_{\max}\in (0, \infty]$, such that
either $R_{\max}=\infty$, or one of the following two
$$
\lim_{r\to R_{max}}\phi(r)=\infty, \qquad \lim_{r\to R_{max}}\phi(r)=0.
$$
In particular, when $\kappa\ge 0$,  $\phi$ is increasing on $[0, R_{\max})$. Then if $R_{\max}<\infty$,
$$
\lim_{r\to R_{max}}\phi(r)=\infty.
$$
\end{lemma}

Before going further, we give the following comparison lemma.

\begin{lemma}\label{lem2-2}
Suppose that $\beta\ge 1$, $A_1\ge A_2$, $\kappa_1\ge\kappa_2$, $A_1+\kappa_1>A_2+\kappa_2$,
  $\phi_1$ and $\phi_2$ solve the following differential inequalities
\begin{align}\label{2-1}
\left\{
\begin{aligned}
&D_v\left(\rho(r)\phi_1'\right)'\ge\kappa_1 \rho(r)\phi_1+A_1\rho(r)\phi_1^{\beta}e^{-\frac{1}{4D_u}r^2},
\\
&D_v\left(\rho(r)\phi_2'\right)'\le\kappa_2 \rho(r)\phi_2+A_2\rho(r)\phi_2^{\beta}e^{-\frac{1}{4D_u}r^2},
\\
\\
&\phi_1(0)=\phi_2(0)=1, \quad \phi_1'(0)=\phi_2'(0)=0.
\end{aligned}\right.
\end{align}
Then $\phi_1>\phi_2$ for any $r>0$. Specifically, when $A_1= A_2$, $\kappa_1=\kappa_2$, if one inequality sign in the above two inequalities is strict, then there are also $\phi_1>\phi_2$ for any $r>0$.
\end{lemma}

{\it\bfseries Proof.}  From L'Hospital principle, we infer that
$$
\phi_1''(0)\ge \frac{A_1+\kappa_1}{ND_v},\quad \phi_2''(0)\le \frac{A_2+\kappa_2}{ND_v}.
$$
Therefore,
$$
\phi_1(r)\ge 1+\frac{A_1+\kappa_1}{2ND_v}r^2+o(r^2), \quad  \phi_2(r)\le 1+\frac{A_2+\kappa_2}{2ND_v}r^2+o(r^2)
$$
for small $r$, which implies that $\phi_1(r)>\phi_2(r)$ for small $r>0$.

We claim that $\phi_1(r)>\phi_2(r)$. Otherwise, there exists $r^*>0$ such that
$$
\phi_1(r)-\phi_2(r)>0, \ \text{for } \ 0<r<r^*;
$$
$$
\phi_1(r^*)-\phi_2(r^*)=0; \quad \phi_1'(r^*)-\phi_2'(r^*)\le 0.
$$
Multiplying both sides of the first inequality of \eqref{2-1} by $\phi_2$, the second inequality by $\phi_1$, integrating them from $0$ to $r^*$, and  then subtracting the second inequality from the first  inequality to get
\begin{align*}
&D_v\left[ \rho(r)\phi_1'\phi_2-\rho(r)\phi_2'\phi_1\right]_{0}^{r^*}
\\
\ge &(\kappa_1-\kappa_2)\int_0^{r^*} \rho(r)\phi_1\phi_2 dr+\int_0^{r^*}\left(A_1\phi_1^{\beta-1}-A_2\phi_2^{\beta-1}\right) \rho(r)e^{-\frac{1}{4D_u}r^2}\phi_1\phi_2dr
\\
> &0,
\end{align*}
that is
$$
0\ge D_v \rho(r^*)\phi_1(r^*)(\phi_1'(r^*)-\phi_2'(r^*))  >0.
$$
It is a contradiction. Therefore, $\phi_1>\phi_2$ for any $r>0$.

In particular, when $A_1= A_2$, $\kappa_1=\kappa_2$, if one inequality sign in the above two inequalities is strict,
the proof is similar to above, we omit it.  . \hfill $\Box$

\medskip

Next, we study the solution of the following linear problem.

\begin{lemma}\label{lem2-3}
Consider
\begin{align}\label{2-2}
\left\{
\begin{aligned}
&D_v\left(r^{N-1}e^{\frac{r^2}{4D_v}}\phi'\right)'=M r^{N-1}e^{\frac{r^2}{4D_v}}\phi,
\\
&\phi(0)=1, \quad \phi'(0)=0.
\end{aligned}\right.
\end{align}
Then $\phi(r)$ exists globally, and there exists a constant $L$ depending on $M$, $D_v$, $N$ such that when $r>1$,
$\phi(r)\le Lr^{4M}$ for $N=1$; and
$\phi(r)\le Lr^{2M}$ for $N\ge 2$.
\end{lemma}

{\it\bfseries Proof.}  By standard fixed point method, it is easy to prove the local existence of \eqref{2-2}. Let $\phi(r)$ be the solution of \eqref{2-2},
and one can also see that $\phi(r)$ is increasing on $r$.  Then
\begin{align*}
D_v r^{N-1}e^{\frac{r^2}{4D_v}}\phi'(r)= &M \int_0^r s^{N-1}e^{\frac{s^2}{4D_v}}\phi(s)ds
\\
\le & M\phi(r)\int_0^r s^{N-1}e^{\frac{s^2}{4D_v}}ds.
\end{align*}
Let $\rho(r)=r^{N-1}e^{\frac{r^2}{4D_v}} $. From the above inequality, we infer that
$$
D_v(\ln\phi(r))'\le M \frac{\int_0^r\rho(s)ds}{\rho(r)}.
$$
Therefore,
\begin{align}\label{eq2-3}
&D_v\ln\phi(r) \le M \int_0^r \frac{\int_0^\tau\rho(s)ds}{\rho(\tau)}d\tau.
\end{align}
In what follows, we calculate $\int_0^r \frac{\int_0^\tau\rho(s)ds}{\rho(\tau)}d\tau$. For this purpose, we first
see that for any large $\tau>1$,
\begin{align}
&\int_0^\tau\rho(s)ds=\int_0^1\rho(s)ds+\frac12\int_1^\tau \frac{1}se^{\frac{s^2}{4D_v}}ds^2\nonumber
\\
&=\int_0^1\rho(s)ds+\left.2D_ve^{\frac{s^2}{4D_v}}\frac{1}s\right|_{1}^{\tau}+2D_v\int_1^\tau \frac{1}{s^2}e^{\frac{s^2}{4D_v}}ds\nonumber
\\
&\le \int_0^1\rho(s)ds+2D_ve^{\frac{\tau^2}{4D_v}}\frac{1}{\tau}+\left.4D_v^2\frac{1}{s^3}e^{\frac{s^2}{4D_v}}\right|_{1}^{\tau}
+12D_v^2\int_1^\tau \frac{1}{s^4}e^{\frac{s^2}{4D_v}}ds\nonumber
\\
&\le \int_0^1\rho(s)ds+2D_ve^{\frac{\tau^2}{4D_v}}\frac{1}{\tau}+4D_v^2\frac{1}{\tau^3}e^{\frac{\tau^2}{4D_v}}
+12D_v^2\frac{1}{\tau^4}e^{\frac{\tau^2}{4D_v}}(\tau-1)\nonumber
\\
\label{eq2-4}
&\le \int_0^1\rho(s)ds+2D_ve^{\frac{\tau^2}{4D_v}}\frac{1}{\tau}+16D_v^2\frac{1}{\tau^3}e^{\frac{\tau^2}{4D_v}}, \ \ \text{when} \  N=1.
\end{align}

\begin{align}\label{eq2-5}
&\int_0^\tau\rho(s)ds=2D_v\left(e^{\frac{\tau^2}{4D_v}}-1 \right), \ \ \text{when} \  N=2.
\end{align}
\begin{align}
&\int_0^\tau\rho(s)ds=\int_0^1\rho(s)ds+\frac12\int_1^\tau s^{N-2}e^{\frac{s^2}{4D_v}}ds^2\nonumber
\\
&=\int_0^1\rho(s)ds+\left.2D_ve^{\frac{s^2}{4D_v}}s^{N-2}\right|_{1}^{\tau}-2(N-2)D_v\int_1^\tau s^{N-3}e^{\frac{s^2}{4D_v}}ds\nonumber
\\
\label{eq2-6}
&\le \int_0^1\rho(s)ds+2D_ve^{\frac{\tau^2}{4D_v}}\tau^{N-2}, \ \ \text{when} \  N\ge 3.
\end{align}
Therefore, substituting \eqref{eq2-4}-\eqref{eq2-6} into \eqref{eq2-3} gives
\begin{align}\label{eq2-7}
\ln\phi(r) \le \frac{M}{D_v} \int_0^1 \frac{\int_0^\tau\rho(s)ds}{\rho(\tau)}d\tau+\frac{M}{D_v}\int_1^r \frac{\int_0^\tau\rho(s)ds}{\rho(\tau)}d\tau
\le 2M\ln r+M_1, \ \ \text{when} \  N=1.
\end{align}
\begin{align}\label{eq2-8}
\ln\phi(r) \le \frac{M}{D_v} \int_0^1 \frac{\int_0^\tau\rho(s)ds}{\rho(\tau)}d\tau+\frac{M}{D_v}\int_1^r \frac{\int_0^\tau\rho(s)ds}{\rho(\tau)}d\tau
\le 2M\ln r+M_N,\ \ \text{when} \  N\ge 2.
\end{align}
This lemma is proved. \hfill $\Box$

In what follows, we show that for any $\kappa\ge 0$,  when $AB^{\alpha-1}$ is appropriately small, the solution must exists globally.

\begin{lemma}\label{lem2-4} Assume $\kappa\ge 0$,  $\alpha+\frac{\chi}{D_u}\ge 1$.
Let $\phi(r)\in C^2[0, R_{\max})$ be the solution of \eqref{1-9}, then when $AB^{\alpha-1}$ is appropriately small, we have $R_{\max}=\infty$,
namely, the solution exists globally. In particular, when $\alpha+\frac{\chi}{D_u}=1$, the solution always exists globally for any $A, B>0$.
\end{lemma}

{\it\bfseries Proof.} Let $\overline\phi$ be the solution of \eqref{2-2} with $M=\kappa+1$.
By Lemma \ref{lem2-2}, we have
\begin{align*}
 \overline\phi\le\left\{
 \begin{aligned}
 & Lr^{4(\kappa+1)}, && \text{for}\ r>1,
 \\
 & L, && \text{for}\ r\le 1.
 \end{aligned}\right.
\end{align*}
Therefore, $e^{-\frac{r^2}{4D_u}} \overline\phi^{\alpha+\frac{\chi}{D_u}-1}$ is bounded.
Take $A B^{\alpha-1}$   appropriately small, such that
$$
A B^{\alpha-1}e^{-\frac{r^2}{4D_u}} \overline\phi^{\alpha+\frac{\chi}{D_u}-1}< 1.
$$
Thus,
$$
D_v\left(\rho\overline\phi'\right)'=(\kappa+1)\rho\overline\phi
> \kappa \rho\overline\phi+A B^{\alpha-1}\rho e^{-\frac{r^2}{4D_u}} \overline\phi^{\alpha+\frac{\chi}{D_u}}.
$$
By Lemma \ref{lem2-2},
$$
\phi(r)\le \overline\phi(r).
$$
Recalling Lemma \ref{lem2-1}, the solution exists globally.

When $\alpha+\frac{\chi}{D_u}=1$, for any $A, B>0$, we set $M=\kappa+A B^{\alpha-1}$, and denote the corresponding solution of \eqref{2-2} by $\phi_M$, from Lemma \ref{lem2-2}, it becomes evident that
$\phi(r)\le \phi_M$, i.e. the solution exists globally.
\hfill $\Box$

 We denote the solution corresponding to the parameter $AB^{\alpha-1}$ by $\phi(r, AB^{\alpha-1})$.
Define
\begin{equation}\label{A}
\mathcal A=\{AB^{\alpha-1}>0, \text{the solution $\phi(r, AB^{\alpha-1})$ exists globally}\},
\end{equation}
and
\begin{equation}\label{B}
\mathcal B=\{AB^{\alpha-1}>0, \text{there exits $R>0$ such that $\lim_{r\to R}\phi(r, AB^{\alpha-1})=\infty$ }\}.
\end{equation}
From Lemma \ref{lem2-2} and Lemma \ref{lem2-4}, we see that when $\kappa\ge 0$, $\mathcal A\cup \mathcal B=\mathbb R^+$ when $\kappa\ge 0$,
and
$$
\sup\mathcal A=\inf\mathcal B.
$$
To show $\mathcal B\not=\emptyset$, we first prove the following lemma.
\begin{lemma}
\label{lem2-5}
Assume that $A>0$, $a>1$, $f, f'>0$ for $r>r_0$, and
$$
f'' \ge Af^{a},\quad \text{ for $r>r_0$}.
$$
Then there exists $R_0>r_0$ such that
$$
\lim_{r\to R_0}f(r)=\infty.
$$
\end{lemma}

{\it\bfseries Proof.} Let $f'=h$. Then
$$
h\frac{dh}{df}>Af^a.
$$
A direct integration yields
$$
h^2\ge \frac{2A}{a+1}\left(f^{a+1}(r)-f^{a+1}(r_0)\right).
$$
That is
$$
f'(r)\ge \sqrt{\frac{2A}{a+1}\left(f^{a+1}(r)-f^{a+1}(r_0)\right)}
$$
Noticing that $f', f''>0$, then there exists $r_1>r_0$ such that when $r>r_1$
$$
f'(r)\ge \sqrt{\frac{A}{a+1}}f^{\frac{a+1}2}(r).
$$
By a direct calculation, we obtain
$$
f^{\frac{1-a}2}(r)\le f^{\frac{1-a}2}(r_1)-\frac{a-1}2\sqrt{\frac{A}{a+1}}(r-r_1).
$$
This lemma is proved since $a>1$.
\hfill $\Box$

Next, we show that $\mathcal B\not=\emptyset$. For this purpose, we let
$$
\theta=\phi \exp\{-\sigma r^2\},
$$
with $\sigma=\frac{1}{4D_u(\alpha-1)+4\chi}$.
Then \eqref{1-8} is transformed into
\begin{align}\label{2-9}\left\{
\begin{aligned}
&D_v\theta''+(4\sigma D_v+\frac12)r\theta'+D_v(N-1)\frac1r\theta'+(4\sigma^2D_v+\sigma)r^2\theta+(2N\sigma D_v-\kappa)\theta=AB^{\alpha-1}\theta^{\alpha+\frac{\chi}{D_u}},
\\
&\theta(0)=1, \theta'(0)=0.
\end{aligned}\right.
\end{align}
It is equivalent to
\begin{align}\label{2-10}\left\{
\begin{aligned}
&D_v(\varrho\theta')'=-(4\sigma^2D_v+\sigma)r^2\varrho\theta-(2N\sigma D_v-\kappa)\varrho\theta
+AB^{\alpha-1}\varrho\theta^{\alpha+\frac{\chi}{D_u}},
\\
&\theta(0)=1, \theta'(0)=0,
\end{aligned}\right.
\end{align}
where $\varrho(r)=r^{N-1}\exp\left\{\left(2\sigma+\frac{1}{4D_v}\right)r^2\right\}$.

In what follows, we shall show that  the  solution $\phi$ may blow up when $AB^{\alpha-1}$ is large.
It is worth noting that this result is applicable to any real number $\kappa$

\begin{lemma}\label{lem2-6} Assume that $\kappa\in \mathbb R$, $\alpha+\frac{\chi}{D_u}>1$.
Let $\phi(r)\in C^2[0, R_{\max})$ be the solution of \eqref{1-8}. Then when $AB^{\alpha-1}$ is appropriately large,   $R_{\max}<\infty$,
that is, the solution $\phi$ blows up at a finite point.
\end{lemma}

{\it\bfseries Proof.} Take $$\underline\theta=1+r^{\lambda}$$ with $\lambda>\max\left\{\frac{2}{\alpha+\frac{\chi}{D_u}-1}, 2\right\}$.
It is straightforward to confirm that $\underline\theta$
serves as a lower solution for \eqref{2-9} if
\begin{align*}
&D_v\lambda(\lambda+N-2)r^{\lambda-2}+(4\sigma D_v+\frac12)\lambda r^{\lambda}+
(4\sigma^2D_v+\sigma)(r^2+r^{\lambda+2})+(2N\sigma D_v-\kappa)(1+r^{\lambda})
\\
&\le AB^{\alpha-1}(1+r^{\lambda})^{\alpha+\frac{\chi}{D_u}}.
\end{align*}
Notably, when $AB^{\alpha-1}$ is sufficiently large, the above inequality is satisfied, thus ensuring that $\underline\theta$ is a lower solution of \eqref{2-9}.
Using Lemma \ref{lem2-2}, it follows that
\begin{equation}\label{2-11}
\theta(r)\ge 1+r^{\lambda}
\end{equation}
when $AB^{\alpha-1}$ is appropriately large. Hence, we have
\begin{align}
D_v(\varrho\theta')'&=\varrho\theta\left(-(4\sigma^2D_v+\sigma)r^2-(2N\sigma D_v-\kappa)
+AB^{\alpha-1}\theta^{\alpha+\frac{\chi}{D_u}-1}\right)\nonumber
\\
\label{2-12}
&\ge \frac12AB^{\alpha-1}\varrho\theta^{\alpha+\frac{\chi}{D_u}},
\end{align}
it implies that $\theta'(r)>0$.
On the other hand, for large $AB^{\alpha-1}$, we also have
\begin{align}\label{2-13}
&D_v\theta''+(4\sigma D_v+\frac12)r\theta'+D_v(N-1)\frac1r\theta'\le \kappa \theta+AB^{\alpha-1}\theta^{\alpha+\frac{\chi}{D_u}}
\le 2AB^{\alpha-1}\theta^{\alpha+\frac{\chi}{D_u}}.
\end{align}
From \eqref{2-13} and $\theta'>0$, we obtain
\begin{align}\label{2-14}
&\left(\frac{D_v}2|\theta'|^2-2AB^{\alpha-1}\frac{1}{\alpha+\frac{\chi}{D_u}+1}\theta^{\alpha+\frac{\chi}{D_u}+1}\right)'\le
-(4\sigma D_v+\frac12)r|\theta'|^2-D_v(N-1)\frac1r|\theta'|^2<0.
\end{align}
Integrating it from $0$ to $r$ yields
$$
\frac{D_v}2|\theta'|^2<2AB^{\alpha-1}\frac{1}{\alpha+\frac{\chi}{D_u}+1}\theta^{\alpha+\frac{\chi}{D_u}+1},
$$
i.e.
\begin{align}\label{2-15}
|\theta'|\le 2\sqrt{AB^{\alpha-1}\frac{1}{\alpha+\frac{\chi}{D_u}+1}}\theta^{\frac{\alpha+1}2+\frac{\chi}{2D_u}}.
\end{align}
From \eqref{2-11} and \eqref{2-15}, we obtain that for large $AB^{\alpha-1}$, when $r>1$,
\begin{align}\label{2-16}
(4\sigma D_v+\frac12)r\theta'+D_v(N-1)\frac1r\theta'<\frac14AB^{\alpha-1}\theta^{\alpha+\frac{\chi}{D_u}},
\end{align}
which together with \eqref{2-12} gives
\begin{align}\label{2-17}
D_v\theta''>\frac14AB^{\alpha-1}\theta^{\alpha+\frac{\chi}{D_u}}, \quad \text{for $r>1$}.
\end{align}
Using Lemma \ref{lem2-5}, we complete the proof. \hfill $\Box$

In what follows, we consider the case $\kappa<0$.

\begin{lemma}\label{lem2-8}
Assume that $\alpha+\frac{\chi}{D_u}\ge 1$, $\kappa<0$.
If $0<AB^{\alpha-1}<2N\sigma D_v-\kappa$, then the solution of \eqref{2-9} is decreasing in $[0, R_{\max})$.
\end{lemma}

{\it\bfseries Proof.} By L'Hospital principle, from \eqref{2-9} we infer that
$$
D_vN\theta''(0)=AB^{\alpha-1}-2N\sigma D_v+\kappa<0.
$$
Then $\theta$ is decreasing in a right neighborhood of $0$.
We claim that $\theta$ is decreasing in $[0, R_{\max})$. Otherwise, there exists a minimal point $r_1$, which is the first minimal point of
$\theta$. Then $\theta(r)$ is decreasing in $(0, r_1)$, and $\theta'(r_1)=0$, $\theta''(r_1)\ge 0$. While from \eqref{2-9}, we derive that
$$
D_v\theta''(r_1)=\theta(r_1)\left(AB^{\alpha-1}\theta^{\alpha+\frac{\chi}{D_u}-1}(r_1)+\kappa-2N\sigma D_v-(4\sigma^2D_v+\sigma)r_1^2\right)<0.
$$
It contradicts to $\theta''(r_1)\ge 0$. Therefore, $\theta$ is decreasing in $[0, R_{\max})$.
This lemma is proved. \hfill $\Box$

Consider the eigenvalue problem
\begin{align}\label{2-18}
\left\{\begin{aligned}
&(\varsigma \varphi' )'+\lambda \varsigma \varphi=0,\qquad r\in(0, R),
\\
&\varphi'(0)=0, \quad \varphi(R)=0.
\end{aligned}\right.
\end{align}
where $ \varsigma(r)=r^{N-1}e^{\delta r^2}$ with $\delta>0$, we denote the first eigenvalue by $\lambda(R)$.

\begin{lemma}\label{lem2-9}
Denote the first eigenvalue of \eqref{2-18} by $\lambda(R)$.

(i)If $R_1>R_2$, then $\lambda(R_1)<\lambda(R_2)$. Also $\lambda(R)\to\infty$ as $R\to 0$.

(ii) $\lambda(\infty)=2N\delta$.
\end{lemma}

The first result (i) is standard (e.g. \cite{CH}). When $R=\infty$,
It is easy to see that the principle eigenfunction $\varphi_1=\exp\left\{-\delta r^2\right\}$,
then the principle eigenvalue $\lambda(\infty)=2N\delta$. Using this, we prove the following lemma.

\begin{lemma}\label{lem2-10}
Assume that $\alpha+\frac{\chi}{D_u}\ge 1$, $-\frac N2\le \kappa<0$. Let $\phi$ be the local solution of \eqref{1-9}. Then $\phi(r)>\exp\left\{- \frac{r^2}{4D_v} \right\}$. In particular, when $\alpha+\frac{\chi}{D_u}=1$, the solution $\phi$  exists globally for any $A, B>0$;
when $\alpha+\frac{\chi}{D_u}>1$, the solution $\phi$  exists globally if  $AB^{\alpha-1}<2N\sigma D_v-\kappa$, and
$$
\lim_{r\to\infty}\theta(r)=0,
$$
where $\theta$ is the solution of \eqref{2-9}, that is $\theta=\phi \exp\left\{-\frac{1}{4D_u(\alpha-1)+4\chi} r^2\right\}$.
\end{lemma}

{\it\bfseries Proof.}Let $\underline\varphi=\exp\left\{- \frac{r^2}{4D_v} \right\}$. Then
we have
$$
D_v\left(\rho\underline\varphi'\right)'=-\frac{N}2\rho\underline\varphi
<\kappa \rho\underline\varphi+A B^{\alpha-1}\rho e^{-\frac{r^2}{4D_u}}\underline\varphi^{\alpha+\frac{\chi}{D_u}}.
$$
By Lemma \ref{lem2-2}, $\phi>\underline\varphi$.

When $\alpha+\frac{\chi}{D_u}=1$, similar to the proof in lemma \ref{lem2-4}, we have $\phi<\phi_M$ for $M>\max\{AB^{\alpha-1}+\kappa, 0\}$.
It implies that $\phi$ exists globally.

When $\alpha+\frac{\chi}{D_u}>1$,  if $AB^{\alpha-1}<2N\sigma D_v-\kappa$, by lemma \ref{lem2-8}, $\theta$ is decreasing in $[0, R_{\max})$.
When $\kappa\ge  -\frac N2$, it is easy to see that $\theta>0$ since $\psi>0$, thus, there exists a constant $\theta^*\ge 0$ such that
$$
\lim_{r\to\infty}\theta(r)=\theta^*.
$$
Therefore, $\lim_{r\to\infty}\theta'(r)=\lim_{r\to\infty}\theta''(r)=0$. Next, we show that $\theta^*=0$. Divide both sides of equation \eqref{2-9} by $r^2$, and then make $r\to\infty$ to obtain
$$
\theta^*=0.
$$
This lemma is proved. \hfill $\Box$

When $\kappa\ge -\frac{N}2$, from Lemma \ref{lem2-4}, Lemma \ref{lem2-6} and Lemma \ref{lem2-10}, we see that the solution exists globally when $AB^{\alpha-1}$ is small,
and the solution blows up at a finite point for large $AB^{\alpha-1}$.
According to Lemma \ref{lem2-2}, $\phi$ is monotonically increasing with respect to $AB^{\alpha-1}$.
Thus, there exists a  critical value for $AB^{\alpha-1}$, we denote it by $A^*$ such that when $AB^{\alpha-1}<A^*$, the solution exists globally, and
when $AB^{\alpha-1}>A^*$, the solution blows up.

In what follows, we show that $A^*\in \mathcal A$.

\begin{lemma}\label{lem2-11}
Assume that  $\alpha+\frac{\chi}{D_u}>1$, $\kappa\ge -\frac N2$.  Let $A^*=\sup\{\mathcal A\}=\inf\{\mathcal B\}$, where $\mathcal A$ and $\mathcal B$ are defined in \eqref{A} and \eqref{B}. Then $A^*\in \mathcal A$.
\end{lemma}

{\it\bfseries Proof.}Suppose to the contrary. Then $\phi(r, A^*)\in \mathcal B$, which implies that there exists $R>0$ such that
$$
\lim_{r\to R^-}\phi(r, A^*)=\infty,\quad \lim_{r\to R^-}\theta(r, A^*)=\infty.
$$
Using the continuous dependency of $\theta$ on $r$ and $AB^{\alpha}$, there exists $\delta>0$ such that \eqref{2-12} holds for
$\theta(r, A^*-\delta)$ when $r\in (R-\delta, R)$. Then similar to the proof of Lemma \ref{lem2-6}, we obtain that $\theta(r, A^*-\delta)$
blows up at a finite point. Then $A^*-\delta\in \mathcal B$, it is a contradiction. \hfill $\Box$

\medskip

With this, the proof of Theorem \ref{thm2-1} is complete.

\section{Existence of self-similar solutions for the case $\kappa< -\frac{N}2$}

In this section,  we consider  the case $\kappa<-\frac N2$, where we observe that when $AB^{\alpha-1}$ is either small or large, the solution does not exist globally. More precisely, we have the following lemma.

\begin{lemma}\label{addlem3-1}
Assume that $\kappa<-\frac N2$, $\alpha+\frac{\chi}{D_u}> 1$. Then when the value of $AB^{\alpha-1}$ is large, the solution blows-up at a finite point;
when the value of $AB^{\alpha-1}$ is small, the solution goes to $0$ at a finite point.
\end{lemma}

{\it\bfseries Proof.}
Recalling Lemma \ref{lem2-6}, we discern that the solution exhibits blow-up behavior at a finite point when the value of $AB^{\alpha-1}$ is large.
Therefore, it suffices to consider the case where the value of $AB^{\alpha-1}$ is small.

Using Lemma \ref{lem2-8}, we see that when $\varsigma=\rho$, the principle eigenvalue
$\lambda(+\infty)=\frac{N}{2D_v}$. Owing to the continuous dependence and monotonicity of $\lambda(R)$ on $R$, we can deduce that for any
$\kappa<-\frac N2$, there exists an $R>0$ such that $-\lambda(R)D_v-\kappa>0$. Let $\varphi_R$ denote the corresponding principle eigenfunction.
Then
\begin{align*}
D_v(\rho \varphi_R' )'=-\lambda(R)D_v \rho \varphi_R\ge \kappa\rho \varphi_R+
+AB^{\alpha-1}e^{\frac{-r^2}{4D_u}}\varphi_R^{\alpha+\frac{\chi}{D_u}}
\end{align*}
if $AB^{\alpha-1}<-\lambda(R)D_v-\kappa$. By Lemma \ref{lem2-2}, $\phi(r)\le \varphi_R(r)$. It implies that $\phi(r)$ reaches $0$ before $R$.
This lemma is proved. \hfill $\Box$

 \medskip

Next, we employ variational methods to select an appropriate value for $AB^{\alpha-1}$, thereby proving the global solvability of Problem \eqref{1-8}.
To achieve this, we begin by introducing several definitions. Let
$$
L_{{\rm rad}, \tilde\rho}^p(\mathbb R^N)=\left\{u; u\ \text{ is spherically symmetric and }\ \int_{\mathbb R^N} \tilde\rho(x) |u|^pdx<\infty\right\},
$$
$$
W_{{\rm rad}, \tilde\rho}^{1,p}=\{u\in L_{{\rm rad}, \tilde\rho}^p(\mathbb R^N); \nabla u\in L_{{\rm rad}, \tilde\rho}^p(\mathbb R^N)\}.
$$
And we define
$$
J(u)=\frac{D_v}2\int_{\mathbb R^N} \tilde\rho(x)|\nabla u|^2dx+\frac{\kappa}2\int_{\mathbb R^N} \tilde\rho(x)|u(x)|^2dx,
$$
$$
H(u)=\frac{1}{q+1}\int_{\mathbb R^N} \tilde\rho(x)\exp\left(-\frac{|x|^2}{4D_u}\right)|u|^{q+1}dx,
$$
with $\tilde\rho(x)=\exp\left(\frac{|x|^2}{4D_v}\right)$.
%Recalling \eqref{1-8}, $\varphi$ is a solution of the equation \eqref{1-8} is equivalent to
Noticing that when $u$ is spherically symmetric, then
$$
\int_{\mathbb R^N} \tilde\rho(x) |u|^pdx=|\partial B_1|\int_0^\infty r^{N-1}\tilde\rho(r) |u(r)|^pdr.
$$
where $|\partial B_1|$ is the measure of the unit sphere surface.
\medskip

We introduce the following  compactness result \cite{KJ}.
\begin{lemma}\label{lem5-1}
Let $\rho=r^{N-1}\tilde\rho(r)$.
Suppose $1<p<q<\infty$ and that
$$
\int_0^1\rho(x)\left(\int_x^1{\rho}^{-\frac1{p-1}}dt\right)^{\frac{r(p-1)}{p}}dx<\infty
$$
for some $r>q$. Then $W_{{\rm rad}, \tilde\rho}^{1,p}$   is compactly embedded in $L_{{\rm rad}, \tilde\rho}^{q}$.
\end{lemma}
Using this lemma, it is not difficult to verify that for $2<p<\frac{2N}{(N-2)_+}$
$$
W_{{\rm rad}, \tilde\rho}^{1,2}\hookrightarrow\hookrightarrow L^{p}_{{\rm rad}, \tilde\rho}.
$$
And we also have
\begin{lemma}\label{lem5-2}
Assume that $q>1$, $u\in W_{{\rm rad}, \tilde\rho}^{1,2}$. For any given $\varepsilon>0$, there exists $C_\varepsilon$ such that
$$
\|u\|_{L^2_{{\rm rad}, \tilde\rho}}^2\le \varepsilon\|\nabla u\|_{L^2_{{\rm rad}, \tilde\rho}}^2+C_\varepsilon\left(\int_0^\infty \rho(r)\exp\left(-\frac{r^2}{4D_u}\right)|u(r)|^{q+1}dr\right)^{\frac{2}{q+1}}.
$$
\end{lemma}

%%%%%%%%%5
%%%%%%%%%%%5
%%%%%%%%%%

{\it\bfseries Proof.} By a direct calculation, it is easy to obtain
\begin{align*}
\int_0^\infty\rho'(r)|u(r)|^2dr=-\rho(0)|u(0)|^2-2\int_0^\infty\rho(r)u(r)u'(r)dr\le \frac{2}{|\partial B_1|}\|u\|_{L^2_{{\rm rad}, \tilde\rho}} \|\nabla u\|_{L^2_{{\rm rad}, \tilde\rho}}.
\end{align*}
Noticing that
$$
\frac{\rho'}{\rho}=\frac{N-1}r+\frac{r}{2D_v}\ge \frac{r}{2D_v},
$$
then for any $R>0$,
\begin{align*}
&\int_0^\infty \rho u^2dr=\int_0^R \rho u^2dr+\int_R^\infty \frac{ \rho}{ \rho'} \rho' u^2dr
\\
&\le \left(\int_0^R \rho(r)\exp\left(-\frac{r^2}{4D_u}\right)|u(r)|^{q+1}dr\right)^{\frac{2}{q+1}}
\left(\int_0^R  \rho(r)\exp\left(\frac{r^2}{2D_u(q-1)}\right)dr\right)^{\frac{q-1}{q+1}}
+\frac{2D_v}{R}\int_R^\infty \rho' u^2dr
\\
&\le \left(\frac{R^N}N\exp\left(\frac{R^2}{4D_v}+\frac{R^2}{2D_u(q-1)}\right)\right)^{\frac{q-1}{q+1}}\left(\int_0^\infty \rho(r)\exp\left(-\frac{r^2}{4D_u}\right)|u(r)|^{q+1}dr\right)^{\frac{2}{q+1}}+\frac{4D_v}{R|\partial B_1|}
\|u\|_{L^2_{{\rm rad}, \tilde\rho}} \|\nabla u\|_{L^2_{{\rm rad}, \tilde\rho}}
\\
&\le C(R)\left(\int_0^\infty \rho(r)\exp\left(-\frac{r^2}{4D_u}\right)|u(r)|^{q+1}dr\right)^{\frac{2}{q+1}}+\frac1{2|\partial B_1|}\|u\|_{L^{2}_{{\rm rad}, \tilde\rho}}^2+
\frac{8D_v^2}{R^2|\partial B_1|}\|\nabla u\|_{L^{2}_{{\rm rad}, \tilde\rho}}^2,
\end{align*}
i.e.
$$
\|u\|_{L^{2}_{{\rm rad}, \tilde\rho}}^2\le 2C(R)|\partial B_1|\left(\int_0^\infty \rho(r)\exp\left(-\frac{r^2}{4D_u}\right)|u(r)|^{q+1}dr\right)^{\frac{2}{q+1}}+
\frac{16D_v^2}{R^2}\|\nabla u\|_{L^{2}_{{\rm rad}, \tilde\rho}}^2.
$$
By the arbitrariness of $R$, we complete the proof. \hfill $\Box$

\begin{lemma}\label{lem5-3}
For any $q<\frac{N+2}{(N-2)_+}$,  $J(u)$ achieves its minimum on the set $\mathcal D=\{u\in W_{{\rm rad}, \tilde\rho}^{1,2}; H(u)=1\}$,
and one of the minimizing functions is non-negative.
\end{lemma}

{\it\bfseries Proof.}
Using Lemma \ref{lem5-2}, and taking $\varepsilon=\frac{D_v}{2|\kappa|}$, it is easy to see that
\begin{align}
J(u)&\ge \frac{D_v}2\|\nabla u\|_{L^{2}_{{\rm rad}, \tilde\rho}}^2+\frac{\kappa}2\left(\frac{D_v}{2|\kappa|}\|\nabla u\|_{L^{2}_{{\rm rad}, \tilde\rho}}^2+C\left(\int_0^\infty  \rho(r)\exp\left(-\frac{r^2}{4D_u}\right)|u(r)|^{q+1}dr\right)^{\frac{2}{q+1}}\right)\nonumber
\\
\label{5-1}
&\ge \frac{D_v}4\|\nabla u\|_{L^{2}_{{\rm rad}, \tilde\rho}}^2-\frac{C|\kappa|}2\left(\int_0^\infty \rho(r)\exp\left(-\frac{r^2}{4D_u}\right)|u(r)|^{q+1}dr\right)^{\frac{2}{q+1}}.
\end{align}
Let
$$
J_0=\inf\{J(u); u\in W_{{\rm rad}, \tilde\rho}^{1,2}, H(u)=1\}.
$$
From \eqref{5-1}, it is easy to see that $J(u)$  is bounded below when $H(u)=1$. Therefore $J_0$ is well-defined.
Let $\{u_n\}_n$  with $u_n\in W_{{\rm rad}, \tilde\rho}^{1,2}$ be a sequence such that $H(u_n)=1$, $J(u_n)\to J_0$ as $n\to\infty$.
We may assume that all the $u_n$ are non-negative since $|u_n|\in W_{{\rm rad}, \tilde\rho}^{1,2}$ and $J(|u_n|)\le J(u_n)$.
It is clear that $J(u_n)$ is bounded. Recalling \eqref{5-1} and Lemma \ref{lem5-2}, then
$$
\|u_n\|_{W_{{\rm rad}, \tilde\rho}^{1,2}}\le C.
$$
Then there exists $u\in W_{{\rm rad}, \tilde\rho}^{1,2}$ such that (by passing to a subsequence, for simplicity, we still denote it by $\{u_n\}$.)
$$
u_n\rightharpoonup u \  \text{in}\  W_{{\rm rad}, \tilde\rho}^{1,2}.
$$
It follows that
\begin{align*}
\|\nabla u\|_{L_{{\rm rad}, \tilde\rho}^{2}}^2\le\liminf_{n\to\infty}\|\nabla u_n\|_{L_{{\rm rad}, \tilde\rho}^{2}}^2,\quad \|u\|_{L_{{\rm rad}, \tilde\rho}^{2}}^2\le \liminf_{n\to\infty}\|u_n\|_{L_{{\rm rad}, \tilde\rho}^{2}}^2.
\end{align*}
Lemma \ref{lem5-1} implies that
$$
u_n\rightarrow u \  \text{in}\  L_{{\rm rad}, \tilde\rho}^p, \ \text{for} \ 2<p<\frac{2N}{(N-2)_+}.
$$
Thus
$$
H(u)=1, \ \text{for} \ q<\frac{N+2}{(N-2)_+}.
$$
In what follows, we show that
\begin{equation}
\label{lim3-2}
\lim_{n\to\infty}\int_0^\infty \rho |u_n-u|^2dr=0.
\end{equation}
Similar to the proof of Lemma \ref{5-2}, for any $p\in \left(2, \frac{2N}{(N-2)_+}\right)$, $R>1$, we derive that
\begin{align*}
&\int_0^\infty \rho |u_n-u|^2dr=\int_0^R \rho |u_n-u|^2dr+\int_R^\infty \frac{ \rho}{ \rho'} \rho' |u_n-u|^2dr
\\
&\le \left(\int_0^R \rho(r)|u_n-u|^pdr\right)^{\frac{2}{p}}
\left(\int_0^R  \rho(r)dr\right)^{\frac{p-2}{p}}
+\frac{2D_v}{R}\int_R^\infty \rho' |u_n-u|^2dr
\\
&\le C_1(R)\left(\int_0^\infty \rho(r)|u_n-u|^pdr\right)^{\frac{2}{p}}+\frac{4D_v}{R|\partial B_1|}
\|u_n-u\|_{L^2_{{\rm rad}, \tilde\rho}} \|\nabla (u_n-u)\|_{L^2_{{\rm rad}, \tilde\rho}}
\\
&\le C_1(R)\left(\int_0^\infty \rho(r)|u_n-u|^pdr\right)^{\frac{2}{p}}+\frac{C_2}{R},
\end{align*}
where $C_2$ is independent of $R$ and $n$. Then for any $\varepsilon>0$, there exists $R_0$ such that when $R\ge R_0$
$$
\frac{C_2}{R_0}<\frac{\varepsilon}2.
$$
Notice that
$$
u_n\rightarrow u \  \text{in} \  L_{{\rm rad}, \tilde\rho}^p,
$$
then there exists $N>0$ such that when $n>N$,
$$
C_1(R_0)\left(\int_0^\infty \rho(r)|u_n-u|^pdr\right)^{\frac{2}{p}}<\frac{\varepsilon}2.
$$
It means that for any $\varepsilon>0$,  when $n>N$,
$$
\int_0^\infty \rho |u_n-u|^2dr<\varepsilon,
$$
and \eqref{lim3-2} is proved. It implies that
$$
\int_0^\infty \rho |u|^2dr=\lim_{n\to\infty}\int_0^\infty \rho |u_n|^2dr.
$$
Summing up, it follows that
\begin{align*}
\|\nabla u\|_{L_{{\rm rad}, \tilde\rho}^{2}}^2\le\liminf_{n\to\infty}\|\nabla u_n\|_{L_{{\rm rad}, \tilde\rho}^{2}}^2,\quad \|u\|_{L_{{\rm rad}, \tilde\rho}^{2}}^2=\lim_{n\to\infty}\|u_n\|_{L_{{\rm rad}, \tilde\rho}^{2}}^2=\liminf_{n\to\infty}\|u_n\|_{L_{{\rm rad}, \tilde\rho}^{2}}^2, \quad H(u)=1.
\end{align*}
Hence,
$$
J_0\le J(u)\le \liminf_{n\to\infty}J(u_n)=\lim_{n\to\infty}J(u_n)=J_0,
$$
which proves this lemma. \hfill $\Box$

Using this lemma, we will prove that the problem \eqref{1-9} admits a global positive solution  for a specific value of $AB^{\alpha-1}$.
\begin{lemma}\label{lem5-4}
Assume that $\kappa<-\frac N2$, $1<q<\frac{N+2}{(N-2)_+}$, then there exists $M^*>0$ such that the following equation
\begin{equation}\label{5-2}
-D_v(\rho{\varphi}')'+\kappa\rho(r) {\varphi}+M^*\rho\exp\left(-\frac{r^2}{4D_u}\right){\varphi}^q=0
\end{equation}
admits a global positive solution  $\varphi^*\in C^2(0,+\infty)$,
and
\begin{equation}\label{5-3}
(\varphi^*)'(r)=\frac{1}{D_v\rho(r)}\int_0^r\left(\kappa\rho(s) {\varphi^*}(s)+M^*\rho(s)\exp\left(-\frac{s^2}{4D_u}\right)|{\varphi^*}(s)|^q\right)ds.
\end{equation}
\end{lemma}

{\it\bfseries Proof.} Let $0\le \varphi^*\in W_{{\rm rad}, \tilde\rho}^{1,2}$ be the minimizing function of $J(u)$ on the set $\mathcal D$.
Using Lagrange multiplier method, there exists a constant $M^*$ such that
$$
DJ(\varphi^*)+M^*DH(\varphi^*)=0,
$$
that is
\begin{equation}\label{5-4}
D_v \int_{\mathbb R^N}\tilde\rho(x) \nabla\varphi^*\nabla v dx+\kappa \int_{\mathbb R^N} \tilde\rho(x) \varphi^*(x)v(x)dx
+M^*\int_{\mathbb R^N} \tilde\rho(x)\exp\left(-\frac{|x|^2}{4D_u}\right)|{\varphi^*}(x)|^{q}v(x)dr=0,
\end{equation}
i.e.
\begin{equation}\label{5-5}
D_v \int_0^\infty \rho(r) {\varphi^*}'(r) v'(r)dr+\kappa \int_0^\infty\rho(r) {\varphi^*}(r)v(r)dr
+M^*\int_0^\infty\rho(r)\exp\left(-\frac{r^2}{4D_u}\right)|{\varphi^*}(r)|^{q}v(r)dr=0,
\end{equation}
 for any $v\in W_{{\rm rad}, \tilde\rho}^{1,2}$ since $\varphi^*$ is spherically symmetric.
Recalling Lemma \ref{lem2-9}, let $\varphi_0=\exp\{-\frac{r^2}{4D_v}\}$ be the principle eigenfunction, $\lambda_0=\frac{N}{2D_v}$ be
the principle eigenvalue, i.e.
$$\left\{
\begin{aligned}
&(\rho\varphi_0')'+\lambda_0\rho\varphi_0=0,
\\
&\varphi_0'(0)=0, \varphi_0(+\infty)=0.
\end{aligned}\right.
$$
Take $v=\varphi_0$ in \eqref{5-5}, it gives
$$
-D_v \int_0^\infty \varphi^*(r) (\rho(r)\varphi_0'(r))'dr+\kappa \int_0^\infty \rho(r) {\varphi^*}(r)\varphi_0(r)dr
+M^*\int_0^\infty\rho(r)\exp\left(-\frac{r^2}{4D_u}\right)|{\varphi^*}(r)|^{q}\varphi_0(r)dr=0,
$$
namely
\begin{equation}\label{5-6}
\left(\frac{N}2+\kappa\right) \int_0^\infty \rho(r) {\varphi^*}(r)\varphi_0(r)dr
+M^*\int_0^\infty \rho(r)\exp\left(-\frac{r^2}{4D_u}\right)|{\varphi^*}(r)|^{q}\varphi_0(r)dr=0,
\end{equation}
Therefore, $M^*>0$ since $\frac{N}2+\kappa<0$, $\varphi_0>0$, $\varphi^*\ge 0$  and $\varphi^*\not\equiv 0$.
From \eqref{5-4}, we see that $\varphi^*$ satisfies \eqref{5-2}
in the sense of  distribution on $(0,\infty)$. Note that $\rho>0$ on $(0,\infty)$, and it is clear that $\varphi^*$ is locally absolutely continuous since
${\varphi^*}'$ is locally integrable on $(0,\infty)$. Hence it is easy to obtain that $\varphi^*\in C^2(0,\infty)$.

Take a function $v$ such that $v=1$ in the neighborhood of $0$ and $v=0$ as $|x|$ large.
Multiplying both sides of the equation \eqref{5-2} by $v$, and then integrating the resulting expression from $x$ to $y$ yields
\begin{align*}
&-D_v\rho(y) {\varphi^*}'(y)v(y)+D_v\rho(x) {\varphi^*}'(x)v(x)+D_v \int_x^y \rho(r) {\varphi^*}'(r) v'(r)dr+\kappa \int_x^y\rho(r) {\varphi^*}(r)v(r)dr
\\
&
+M^*\int_x^y\rho(r)\exp\left(-\frac{r^2}{4D_u}\right)|{\varphi^*}(r)|^{q}v(r)dr=0.
\end{align*}
Letting $y\to\infty$, $x\to 0$ yields
\begin{align*}
&D_v\rho(0) {\varphi^*}'(0)+D_v \int_0^{\infty} \rho(r) {\varphi^*}'(r) v'(r)dr+\kappa\int_0^{\infty}\rho(r) {\varphi^*}(r)v(r)dr
\\
&
+M^*\int_0^{\infty}\rho(r)\exp\left(-\frac{r^2}{4D_u}\right)|{\varphi^*}(r)|^{q}v(r)dr=0.
\end{align*}
From \eqref{5-5}, it is easy to get that
$$
\rho(0) {\varphi^*}'(0)=0.
$$
Hence, \eqref{5-3} holds.
Since $\varphi^*\ge 0$  and $\varphi^*\not\equiv 0$, if there exists an $r_0>0$ such that $\varphi^*(r_0)=0$, it can be easily deduced from the  equation \eqref{5-2} that  ${\varphi^*}'(r_0)<0$, and this will cause the sign of $\varphi^*$ to change, which contradicts the non-negativity of $\varphi^*$. Consequently, $\varphi^*(r)>0$ for all $r\in(0,\infty)$. The proof of this lemma is complete. \hfill $\Box$

In the aforementioned lemma, we have proven the global existence of positive solutions for equation \eqref{5-2} in the interval $(0, +\infty)$. Next, we will demonstrate the continuity and differentiability of the solution at point $0$.

\begin{lemma}
\label{lem5-5}
Assume that $\kappa<-\frac N2$, $1<q<\frac{N+2}{(N-2)_+}$.
Let $\varphi^*\in C^2(0, +\infty)$  be the global positive solution of \eqref{5-2}. Then
$\varphi^*\in C^2[0, +\infty)$, and $\varphi^*(0)>0$, $(\varphi^*)'(0)=0$.
\end{lemma}

{\it\bfseries Proof.} We first show that $\varphi^*\in C[0, +\infty)$. To achieve this, it suffices to prove that $(\varphi^*)'\in L^1\left(0, \frac12\right)$.
Note that $(\varphi^*)'\in L^{2}_{\rho}$ with $\rho=r^{N-1}\exp(\frac{|r^2|}{4D_v})$, then it is clearly true when $N=1$.

Next, we consider the cases $N\ge 2$.
Recalling \eqref{5-3}, we see that for any $r\in(0, \frac12)$,
\begin{align}
|(\varphi^*)'(r)| &\le \frac{1}{D_v\rho(r)}\int_0^r\left(|\kappa|\rho(s) {\varphi^*}(s)+M^*\rho(s)\exp\left(-\frac{s^2}{4D_u}\right)|{\varphi^*}(s)|^q\right)ds\nonumber
\\
&\le \frac{1}{D_v\rho(r)}|\kappa|\left(\int_0^r\rho|\varphi^*|^{p}ds\right)^{\frac1{p}}
\left(\int_0^r\rho(s)ds\right)^{\frac{p-1}{p}}+M^*\left(\int_0^r\rho(s)|{\varphi^*}(s)|^{p}ds\right)^{\frac{q}p}
\left(\int_0^r\rho(s)ds\right)^{\frac{p-q}p}\nonumber
\\
&\le C\left(r^{1-\frac Np}+r^{1-\frac{q N}p}\right)\nonumber
\\
\label{5-7}
&\le \tilde Cr^{1-\frac{q N}p}
\end{align}
for any $q<p<\frac{2N}{(N-2)_+}$.

i) When $N=2$, take $p=qN$, from \eqref{5-7}, it is easy to see that
\begin{align}\label{5-8}
|(\varphi^*)'(r)| \le C, \ \text{for} \ r\in \left(0, \frac12\right),
\end{align}
Then $(\varphi^*)'\in L^1\left(0, \frac12\right)$.

ii) When $N\ge 3$, noticing that $q<\frac{N+2}{N-2}$, it means that $q+1<\frac{2N}{N-2}$. Take $p=q+1$, then
\begin{align}\label{5-9}
|(\varphi^*)'(r)|\le \tilde Cr^{1-\frac{q N}{q+1}}
\end{align}
noticing that when $q<\frac{2}{N-2}$, $1-\frac{q N}{q+1}>-1$, then $(\varphi^*)'\in L^1\left(0, \frac12\right)$.
In what follows, it suffices to consider the case $1-\frac{q N}{q+1}\le -1$, i.e. $\frac{2}{N-2}\le q<\frac{N+2}{N-2}$.

Noticing that if
\begin{align}\label{5-10}
|(\varphi^*)'(r)|\le \tilde Cr^{-1-\gamma}, \ \text{for }\ r\in \left(0, \frac12\right),
\end{align}
for $\gamma\ge 0$
then
\begin{align}\label{5-11}
\varphi^*(r)\le \varphi^*(1)+\int_r^1|(\varphi^*)'(s)|ds\le
\left\{\begin{aligned}
&Cr^{-\gamma}, &&\text{for $\gamma>0$},
\\
&C|\ln r|, && \text{for $\gamma=0$}.
\end{aligned}\right.
\end{align}
Using \eqref{5-3}, when $\gamma=0$, we conclude that
\begin{align}
|(\varphi^*)'(r)| &\le \frac{C}{\rho(r)}\int_0^r\left(\rho(s)|\ln s|^q\right)ds\nonumber
\\
&\le  \frac{C}{\rho(r)}\int_0^r s^{N-1-\eta}ds\nonumber
\\
\label{5-12}
&\le Cr^{1-\eta}
\end{align}
for any $\eta\in (0, 1)$, it implies that $(\varphi^*)'\in L^1\left(0, \frac12\right)$.

When $0<\gamma<\frac{N}{q}$,
\begin{align}
|(\varphi^*)'(r)| &\le \frac{C}{\rho(r)}\int_0^r s^{N-1-\gamma q}ds\nonumber
\\
\label{5-13}
&\le Cr^{1-\gamma q}=Cr^{-1-(\gamma q-2)}.
\end{align}
Noting that $1-\gamma q>-1-\gamma$  is equivalent to $\gamma<\frac{2}{q-1}$. In summary, from \eqref{5-10}, \eqref{5-11} and \eqref{5-13}, we conclude that if
\begin{equation}\label{5-14}
\gamma<\min\left\{\frac{N}{q}, \frac{2}{q-1}\right\},
\end{equation}
then \eqref{5-13} can be deduced from \eqref{5-10}, such that $\gamma q-2<\gamma$.
We define the sequence $\{\gamma_n\}_n$ with $\gamma_{n+1}=\gamma_n q-2$, $\gamma_1<\min\left\{\frac{N}{q}, \frac{2}{q-1}\right\}$.
It is clear that $\gamma_n$ is decreasing, and there exists $K>0$ such that $\gamma_K<0$, it implies that
$(\varphi^*)'\in L^1\left(0, \frac12\right)$.

Therefore, when $N\ge 3$, if \eqref{5-10} holds, then $(\varphi^*)'\in L^1\left(0, \frac12\right)$ is ensured by \eqref{5-14}.
Recalling \eqref{5-9}, we see that $\gamma=\frac{q N}{q+1}-2$.
Hence, it suffices to verify that
$$
\frac{q N}{q+1}-2<\min\left\{\frac{N}{q}, \frac{2}{q-1}\right\}.
$$
It is easy to see that $\frac{q N}{q+1}-2=N-2-\frac{N}{q+1}$, by a direct calculation, we see that
$$
\frac{q N}{q+1}-2<\frac{N}{q}\Leftrightarrow \frac{N-2}{N}<\frac1{q+1}+\frac1q,
$$
$$
\frac{q N}{q+1}-2<\frac{2}{q-1}\Leftrightarrow N-2<\frac{2}{q-1}+\frac{N}{q+1}.
$$
It is clear that the above two inequalities can be ensured by $q<\frac{N+2}{N-2}$. Thus, $\varphi^*\in C[0, +\infty)$ is proved.

Recalling \eqref{5-3}, and noticing that $\varphi^*\in C[0, +\infty)$, then
$$
|(\varphi^*)'(r)|\le Cr,\quad \text{for $r<1$}.
$$
Hence $(\varphi^*)'(r)\in C[0,\infty)$ with $(\varphi^*)'(0)=0$. Using \eqref{5-3}, we further have $(\varphi^*)''(r)\in C[0,\infty)$.

At last, we show that $\varphi^*(0)>0$.
Recalling Lemma \ref{lem2-8}, we see that when $\varsigma=\rho$, the principle eigenvalue
$\lambda(+\infty)=\frac{N}{2D_v}$. From the continuous dependence and monotonicity of $\lambda(R)$ on $R$, it can be concluded that
for any $\kappa<-\frac N2$, there exists $R>0$ such that $-\lambda(R)D_v-\kappa>0$. Let $\varphi_R$ be the corresponding principle eigenfunction.
Then
\begin{align*}
D_v(\rho \varphi_R' )'=-\lambda(R)D_v \rho \varphi_R\ge \kappa\rho \varphi_R+
M^*\rho\exp\left(-\frac{r^2}{4D_u}\right){\varphi_R}^q
\end{align*}
when $\varphi_R(0)<\left|\frac{-\lambda(R)D_v-\kappa}{M^*}\right|^{\frac1{q-1}}$. Therefore, if $\varphi^*(0^+)\le \varphi_R(0)$, by Lemma \ref{lem2-2}, $\varphi^*(r)\le \varphi_R(r)$, and $\varphi^*(r)$ reaches $0$ before $r=R$. By a direct analysis, it is also easy to obtain that
$\varphi^*(r)$ changes sign after $R$, it leads to a contradiction. Hence, $\varphi^*(0)>0$.
 We complete the proof. \hfill $\Box$

 \medskip

 {\it\bfseries Proof of Theorem \ref{thm2-2}. }  By Lemma \ref{lem5-4} and Lemma \ref{lem5-5}, when $\kappa<-\frac N2$, $1<q<\frac{N+2}{(N-2)_+}$,
the equation \eqref{5-2} admits  a global positive solution $\varphi^*\in C^2[0, +\infty)$ such that $\varphi^*(0)>0$, $(\varphi^*)'(0)=0$.
Let $\tilde\phi^*(r)=\frac{\varphi^*(r)}{\varphi^*(0)}$. Then $\tilde\phi^*$ satisfies \eqref{1-9} with $AB^{\alpha-1}=M^*|\varphi^*(0)|^{\alpha+\frac{\chi}{D_u}-1}$ and the global solvability for $1<\alpha+\frac{\chi}{D_u}<\frac{N+2}{(N-2)_+}$ is thus proven. On the other hand, since $\tilde\phi^*\in C^2[0,\infty)\cap W_{\rho}^{1,2}(0,\infty)$, then
$$
\lim_{r\to\infty}\rho(r)|\phi(r)|^2=0.
$$
Hence, $\tilde\phi^*(r)=o(\rho^{-\frac12})$, as $r\to\infty$. We complete the proof.
\hfill $\Box$

\section{Growth or decay rate estimation, and  classification of self-similar singular solutions}

In Section 2, we study the global existence of solutions to Problem  \eqref{1-8} (or \eqref{1-9}) for  the case $\kappa\ge-\frac N2$, proving that when $AB^{\alpha-1}$  is small, Problem \eqref{1-8} has a global solution; whereas when $AB^{\alpha-1}$  is large, the solution blows up in finite time.
In Section 3,  we consider  the case $\kappa<-\frac N2$, where we observe that when $AB^{\alpha-1}$ is either small or large, the solution does not exist globally. We employ variational methods to select an appropriate value for $AB^{\alpha-1}$, thereby proving the global solvability of Problem \eqref{1-8}.

In this section, our primary focus is on estimating the growth  or decay rate of the global solution at infinity.
We first consider the case $\kappa\ge-\frac N2$, from the proof of Lemma \ref{lem2-4},  we speculate that $\phi$ may have an algebraic growth rate for small $AB^{\alpha-1}$. Therefore, we make the transformation as follows. Let
$$
\phi(r)=r^\eta g_{\eta}(r).
$$
Substituting it into \eqref{1-8} yields
\begin{equation}
\label{4-1}
D_vg_{\eta}''+\left(\frac{r}2+\frac{D_v(2\eta+N-1)}r\right)g_{\eta}'=\left(\kappa-\frac{\eta}2-\frac{\eta D_v(\eta+N-2)}{r^2}\right)g_{\eta}+AB^{\alpha-1}r^{\eta(\alpha+\frac{\chi}{D_u}-1)}g_{\eta}^{\alpha+\frac{\chi}{D_u}}e^{-\frac{r^2}{4D_u}}.
\end{equation}
From Section 2, we see that when $AB^{\alpha-1}$ is appropriately small, the solution $g_{\eta}$ with $g_{\eta}(r)>0$ exists globally.
In order to estimate the growth or decay rate of $g_{\eta}$, we first prove the following lemma.

\begin{lemma}\label{lem4-1}
Assume that $\kappa\ge-\frac N2$. Let $\phi\in C^2[0,\infty)$ be the global solution of \eqref{1-8}.
When $\alpha+\frac{\chi}{D_u}>1$, and $AB^{\alpha-1}$ is appropriately small, for any $\eta>2\kappa$,  the function $g_{\eta}(r)$ is decreasing for large $r$, and
$$
\lim_{r\to\infty} g_{\eta}(r)=0.
$$
When $\alpha+\frac{\chi}{D_u}=1$, the above limit relation holds true
for any positive $A$ and $B$.
\end{lemma}

{\it\bfseries Proof.} Since  $\eta>2\kappa$, we claim that there exists a sufficient large positive constant $r_0$ with
$$
\kappa-\frac{\eta}2-\frac{\eta D_v(\eta+N-2)}{r_0^2}<\frac12\left(\kappa-\frac{\eta}2\right)<0
$$
such that
\begin{equation}\label{4-2}
g_{\eta}'(r)<0,\quad \text{for}\quad r>r_0.
\end{equation}

Otherwise, one of the following two situations must occur:

\medskip

(i) $g_{\eta}$ oscillates when $r$ is sufficiently large.

(ii) There exists $r_1>r_0$ such that $g_{\eta}$ is increasing for $r>r_1$.

\medskip

If (i) holds. Then there exists a sequence of minimum points  $\{r_n\}$ of $g_{\eta}$, such that $r_n>r_0$
and $r_n\to \infty$ as $n\to\infty$.
Then $g_{\eta}'(r_n)=0$, $g_{\eta}''(r_n)\ge 0$.
While by \eqref{4-1}, we also have
\begin{align*}
D_vg_{\eta}''(r_n)&=\left(\kappa-\frac{\eta}2-\frac{\eta D_v(\eta+N-2)}{r^2}\right)g_{\eta}(r_n)+AB^{\alpha-1}r^{\eta(\alpha+\frac{\chi}{D_u}-1)}g_{\eta}^{\alpha+\frac{\chi}{D_u}}(r_n)e^{-\frac{r^2}{4D_u}}
\\
&\le \left(\frac12\left(\kappa-\frac{\eta}2\right)+AB^{\alpha-1}r_n^{\eta(\alpha+\frac{\chi}{D_u}-1)}
g_{\eta}^{\alpha+\frac{\chi}{D_u}-1}(r_n)e^{-\frac{r^2}{4D_u}}\right)g_{\eta}(r_n)
\end{align*}
From the proof of Lemma \ref{lem2-4}, we see that when $\kappa\ge 0$, $\phi$ takes a function with algebraic growth rate as the upper bound. When $\kappa<0$,  this conclusion also holds according to the comparison lemma \ref{lem2-2}. Hence, when $r_n$ is sufficiently large, we have
$$
AB^{\alpha-1}r_n^{\eta(\alpha+\frac{\chi}{D_u}-1)}g_{\eta}^{\alpha+\frac{\chi}{D_u}-1}(r_n)e^{-\frac{r^2}{4D_u}}<\frac18\left(\eta-2\kappa\right).
$$
Combining the above two inequalities, when $n$ is sufficiently large, we have
\begin{equation*}
D_vg_{\eta}''(r_n)<\frac14\left(\kappa-\frac{\eta}2\right)g_{\eta}(r_n)<0.
\end{equation*}
It contradicts to $g_{\eta}''(r_n)\ge 0$.

If (ii) holds, that is $g_{\eta}'(r)\ge 0$ for $r>r_1$.
From \eqref{4-1}, there exists $r_2>r_1$ such that when $r>r_2$
\begin{align*}
D_vg_{\eta}''=&-\left(\frac{r}2+\frac{D_v(2\eta+N-1)}r\right)g_{\eta}'+\left(\kappa-\frac{\eta}2-\frac{\eta D_v(\eta+N-2)}{r^2}\right)g_{\eta}+AB^{\alpha-1}r^{\eta(\alpha+\frac{\chi}{D_u}-1)}g_{\eta}^{\alpha+\frac{\chi}{D_u}}e^{-\frac{r^2}{4D_u}}
\\
\le &\frac12\left(\kappa-\frac{\eta}2\right)g_{\eta}+AB^{\alpha-1}r^{\eta(\alpha+\frac{\chi}{D_u}-1)}g_{\eta}^{\alpha+\frac{\chi}{D_u}}e^{-\frac{r^2}{4D_u}}
\\
\le & \frac14\left(\kappa-\frac{\eta}2\right)g_{\eta}
\end{align*}
since $g_{\eta}$  has an upper bound function with algebraic growth rate.
Integrating the above inequality from $r_1$ to $r$, we obtain
$$
D_vg_{\eta}'(r)<D_vg_{\eta}'(r_2)+\frac14\left(\kappa-\frac{\eta}2\right)\int_{r_2}^rg_{\eta}(s)ds<D_vg_{\eta}'(r_2)+
\frac14\left(\kappa-\frac{\eta}2\right)g_{\eta}(r_2)(r-r_2)\to-\infty
$$
as $r\to\infty$. It contradicts to $g_{\eta}'\ge 0$.
Summing up, we conclude \eqref{4-2}. It implies that $g_{\eta}(r)$ with $g_{\eta}(r)>0$ is decreasing. Therefore $\lim_{r\to\infty}g_{\eta}(r)$ exists.
Recalling \eqref{4-1}, and letting $r\to\infty$, we obtain
$$
\lim_{r\to\infty} g_{\eta}(r)=0.
$$
When $\alpha+\frac{\chi}{D_u}=1$, it is easy to see that
$$
D_vg_{\eta}''+\left(\frac{r}2+\frac{D_v(2\eta+N-1)}r\right)g_{\eta}'<\left(\kappa-\frac{\eta}2\right)g_{\eta},
$$
similar to the proof above, this lemma is proved. \hfill $\Box$

From Lemma \ref{lem4-1}, we can observe that for any positive $\varepsilon$,  $\phi(r)<Cr^{2\kappa+\varepsilon}$ for large $r$.
Next, we turn our attention to estimating the lower bound of $\phi(r)$. Initially, we consider the case where  $\kappa\ge 0$.

\begin{lemma}\label{lem4-2}
Assume that $\kappa\ge 0$, $\alpha+\frac{\chi}{D_u}\ge 1$.
Let $\phi(r)\in C^2[0,\infty)$ be the global positive solution of \eqref{1-8}. Then when $AB^{\alpha-1}$ is appropriately small, there exists a positive constant $C^*$ such that
$$
\phi(r)\ge C^* r^{2\kappa},
$$
for large $r$.
In particular, when $\alpha+\frac{\chi}{D_u}=1$, the above estimate holds for any $A, B>0$.
\end{lemma}

{\it\bfseries Proof.} By performing a direct integration from $0$ to $r$ of  the equation \eqref{1-9}, it is straightforward to deduce that
$$
D_v \rho(r)\phi'(r)=\int_0^r\left(\kappa \rho\phi+A B^{\alpha-1}\rho e^{-\frac{r^2}{4D_u}} \phi^{\alpha+\frac{\chi}{D_u}}\right) ds>0
$$
due to the condition  $\kappa\ge 0$. Consequently, $\phi'(r)>0$ for $r>0$.
According to Lemma \ref{lem4-1},  when $AB^{\alpha-1}$ is appropriately small,
the speed of $\phi$ at infinity grows at most algebraically, Hence, there exists an $M>0$ such that
\begin{align*}
0<D_v\left(\rho\phi'\right)'=\kappa \rho\phi+A B^{\alpha-1}\rho e^{-\frac{r^2}{4D_u}} \phi^{\alpha+\frac{\chi}{D_u}}<M\rho\phi.
\end{align*}
It is also noteworthy that when $\alpha+\frac{\chi}{D_u}=1$,
the above inequality holds true for any $A, B$.
Integrating the above inequality from $0$ to $r$ gives
$$
0<D_v\rho\phi'(r)<M\int_0^r\rho\phi ds<M\phi(r)\int_0^r\rho ds,
$$
that is
\begin{align}\label{eq4-3}
0<(\ln \phi(r))'<\frac{M}{D_v}\frac{1}{\rho(r)}\int_0^r\rho ds.
\end{align}
Observing that
$$
\lim_{r\to\infty}\frac{M}{D_v}\frac{\int_0^r\rho ds}{r^{-1}\rho(r)}=\lim_{r\to\infty}\frac{M}{D_v}\frac{\rho(r)}{\frac{\rho}{2D_v}+(N-2)r^{-2}\rho}=2M,
$$
we conclude that there exists an $r_0>1$ such that
$$
\frac{M}{D_v}\frac{1}{\rho(r)}\int_0^r\rho ds<\frac{3M}{r}, \ \text{when $r>r_0$}.
$$
Combining this with \eqref{eq4-3}, we obtain
\begin{align}\label{eq4-4}
0<(\ln \phi(r))'<\frac{3M}{r}, \ \text{when $r>r_0$}.
\end{align}
Next,  we define $h(r)=\ln{\phi(r)}$. Then
\begin{align}\label{eq4-5}
\left\{
\begin{aligned}
&\kappa -\frac{1}2rh'(r)=D_v\left(h''(r)+\frac{N-1}{r}h'(r)+|h'(r)|^2\right)-AB^{\alpha-1}
e^{\left(\alpha+\frac{\chi}{D_u}-1\right)h}e^{-\frac{1}{4D_u}r^2},
\\
&h(0)=0, \quad h'(0)=0.
\end{aligned}\right.
\end{align}
From \eqref{eq4-3} and \eqref{eq4-4}, we derive that
\begin{align*}
D_v(\rho(r)h')'=-D_v|h'|^2\rho(r)+\kappa\rho(r)+AB^{\alpha-1}\rho(r)
e^{\left(\alpha+\frac{\chi}{D_u}-1\right)h}e^{-\frac{1}{4D_u}r^2}
>\kappa\rho(r)-\frac{9D_vM^2}{r^2}\rho(r), \ \text{when $r>r_0$}.
\end{align*}
Given that $h(0)=\ln\phi(0)=0$, $h'(r)>0$,  it follows that $h(r)>0$.
Integrating the above inequality from $r_0$ to $r$ gives
 \begin{align}
h(r)&>h(r_0)+\int_{r_0}^r\frac{\rho(r_0)h'(r_0)}{\rho(s)}ds+\frac{\kappa}{D_v}\int_{r_0}^r\frac{1}{\rho(\tau)}\int_{r_0}^{\tau}\rho(s)dsd\tau
-9M^2\int_{r_0}^r\frac{1}{\rho(\tau)}\int_{r_0}^{\tau}s^{-2}\rho(s)dsd\tau\nonumber
\\
\label{eq4-6}
&>\frac{\kappa}{D_v}\int_{r_0}^r\frac{1}{\rho(\tau)}\int_{r_0}^{\tau}\rho(s)dsd\tau
-9M^2\int_{r_0}^r\frac{1}{\rho(\tau)}\int_{r_0}^{\tau}s^{-2}\rho(s)dsd\tau.
\end{align}
Next, we focus our attention on calculating
 $\int_{r_0}^r \frac{\int_{r_0}^\tau\rho(s)ds}{\rho(\tau)}d\tau$. Initially, we observe that
\begin{align}
\int_{r_0}^\tau\rho(s)ds&=\int_{r_0}^{\tau}s^{N-1}e^{\frac{s^2}{4D_v}}ds\nonumber
\\
&=\left.2D_ve^{\frac{s^2}{4D_v}}s^{N-2}\right|_{r_0}^{\tau}-2D_v(N-2)\int_{r_0}^\tau s^{N-3}e^{\frac{s^2}{4D_v}}ds\nonumber
\\
&= 2D_ve^{\frac{\tau^2}{4D_v}}\tau^{N-2}-2D_ve^{\frac{r_0^2}{4D_v}}{r_0}^{N-2}-2D_v(N-2)\int_{r_0}^\tau s^{N-3}e^{\frac{s^2}{4D_v}}ds.
\label{eq4-7}
\end{align}
Therefore,
\begin{align}\label{eq4-8}
&\int_{r_0}^\tau\rho(s)ds\ge 2D_ve^{\frac{\tau^2}{4D_v}}\tau^{N-2}-2D_ve^{\frac{r_0^2}{4D_v}}{r_0}^{N-2}, \quad \text{for $N=1, 2$}.
\end{align}
When $N=3$, noticing that $r_0>1$, then from \eqref{eq4-7} we infer that
\begin{align}
\int_{r_0}^\tau\rho(s)ds&> 2D_ve^{\frac{\tau^2}{4D_v}}\tau-2D_ve^{\frac{r_0^2}{4D_v}}{r_0}-2D_v\int_{r_0}^\tau s e^{\frac{s^2}{4D_v}}ds\nonumber
\\
&\ge 2D_ve^{\frac{\tau^2}{4D_v}}\tau -2D_ve^{\frac{r_0^2}{4D_v}}{r_0}-4D_v^2e^{\frac{\tau^2}{4D_v}}+4D_v^2e^{\frac{{r_0}^2}{4D_v}}, \quad \text{for $N=3$}.
\label{4-9}
\end{align}
For higher-dimensional cases $N\ge 4$, leveraging \eqref{eq4-7},  we further obtain
\begin{align}
\int_{r_0}^\tau\rho(s)ds&= 2D_ve^{\frac{\tau^2}{4D_v}}\tau^{N-2}-2D_ve^{\frac{r_0^2}{4D_v}}{r_0}^{N-2}-2D_v(N-2)\int_{r_0}^\tau s^{N-3}e^{\frac{s^2}{4D_v}}ds\nonumber
\\
&=2D_ve^{\frac{\tau^2}{4D_v}}\tau^{N-2}-2D_ve^{\frac{r_0^2}{4D_v}}{r_0}^{N-2}-4D_v^2(N-2)e^{\frac{\tau^2}{4D_v}}\tau^{N-4}+4D_v^2(N-2)
e^{\frac{{r_0}^2}{4D_v}}{r_0}^{N-4}\nonumber
\\
&\qquad+4D_v^2(N-2)(N-4)
\int_{r_0}^\tau s^{N-5}e^{\frac{s^2}{4D_v}}ds\nonumber
\\
&\ge 2D_ve^{\frac{\tau^2}{4D_v}}\tau^{N-2}-2D_ve^{\frac{r_0^2}{4D_v}}{r_0}^{N-2}-4D_v^2(N-2)e^{\frac{\tau^2}{4D_v}}\tau^{N-4}.
\label{4-10}
\end{align}
Combining \eqref{eq4-8}, \eqref{4-9} and \eqref{4-10}, we derive that
\begin{align}\label{eq4-11}
\frac{\kappa}{D_v}\int_{r_0}^r\frac{1}{\rho(\tau)}\int_{r_0}^{\tau}\rho(s)dsd\tau\ge  2\kappa\ln r-K_1,
\end{align}
Similarly, we also establish the following boundedness estimate: there exists a constant $K_2>0$ such that
\begin{align}\label{eq4-12}
9M^2\int_{r_0}^r\frac{1}{\rho(\tau)}\int_{r_0}^{\tau}s^{-2}\rho(s)dsd\tau<K_2.
\end{align}
Substituting \eqref{eq4-11} and \eqref{eq4-12} into \eqref{eq4-6}, we obtain
$$
h(r)>2\kappa\ln r-K_1-K_2.
$$
where $K_1, K_2$ are  constants. This completes the proof of the lemma.
 \hfill $\Box$

\medskip

Next, we consider the case $-\frac{N}2\le \kappa<0$. To derive the estimate presented in Lemma \ref{lem4-2},
this situation proves to be significantly more intricate than the case where $\kappa\ge 0$.
To accomplish this objective, we initially establish the following lower bound decay rate estimate for
 $\phi$.

\begin{lemma}\label{addlem4-3}
Assume that  $-\frac{N}2\le \kappa<0$, and $\alpha+\frac{\chi}{D_u}\ge 1$.
Let $\phi(r)\in C^2[0,\infty)$ be the global positive solution of \eqref{1-8}. Then for any $\eta<2\kappa$,
there exists a positive constant $C^*$ such that
$$
\phi(r)>Cr^{\eta}
$$
for large $r$.
\end{lemma}

{\it\bfseries Proof.} Let $\phi(r)={\bar g}(r)\exp\{-\frac{r^2}{4D_v}\}$. Then
$$
D_v (r^{N-1}{\bar g}')'-\frac{1}{2}r^N{\bar g}'=\left(\frac{N}{2}+\kappa\right)r^{N-1}{\bar g}+A B^{\alpha-1}\rho e^{-\frac{r^2}{4D_u}} \phi^{\alpha+\frac{\chi}{D_u}},
$$
therefore,
\begin{align}\label{eq4-13}\left\{\begin{aligned}
&D_v \left(\bar\rho(r){\bar g}'\right)'=\left(\frac{N}{2}+\kappa\right)\bar\rho(r){\bar g}
+A B^{\alpha-1}\bar\rho(r){\bar g}^{\alpha+\frac{\chi}{D_u}}
\exp\left\{-\left(\alpha+\frac{\chi}{D_u}-1\right)\frac{r^2}{4D_v}-\frac{r^2}{4D_u}\right\},
\\
&{\bar g}(0)=1, {\bar g}'(0)=0,
\end{aligned}\right.\end{align}
where $\bar\rho(r)=r^{N-1}\exp\left\{-\frac{r^2}{4D_v}\right\}$. Noticing that $\frac{N}{2}+\kappa\ge 0$, from the above equation, we derive that
$$
\left(\bar\rho(r){\bar g}'\right)'>0,
$$
which implies that ${\bar g}'(r)>0$ for any $r>0$. By  direct integration, we obtain that
$$
\bar\rho(r){\bar g}'(r)>{\bar g}'(1)\bar\rho(1)>0, \ \text{when $r>1$},
$$
that is
$$
{\bar g}'(r)>A_0r^{1-N}\exp\left\{\frac{r^2}{4D_v}\right\},
$$
where $A_0={\bar g}'(1)\bar\rho(1)>0$.
Therefore,
\begin{align*}
{\bar g}(r)>A_0\int_1^rs^{1-N}\exp\left\{\frac{s^2}{4D_v}\right\}ds
>2D_vA_0r^{-N}\left(\exp\left\{\frac{r^2}{4D_v}\right\}-\exp\left\{\frac{1}{4D_v}\right\}\right).
\end{align*}
It means that when $r>1$,
$$
\phi(r)>2D_vA_0r^{-N}\left(1-\exp\left\{\frac{1-r^2}{4D_v}\right\}\right).
$$
Hence, there exists a positive constant $C^*$, such that
\begin{align}\label{eq4-14}
\phi(r)\ge C^* r^{-N}, \ \text{ when $r>2$}.
\end{align}

Recalling \eqref{4-1} with $\phi(r)=r^\eta g_{\eta}(r)$, then for any $\eta<2\kappa$, there exists $r_1>2$ such that when $r>r_1$,
$$
D_vg_{\eta}''+\left(\frac{r}2+\frac{D_v(2\eta+N-1)}r\right)g_{\eta}'\ge \frac12\left(\kappa-\frac{\eta}2\right)g_{\eta}.
$$
It is equivalent to
\begin{equation}\label{eq4-15}
D_v\left(\tilde\rho g_{\eta}'\right)'\ge \frac12\left(\kappa-\frac{\eta}2\right)\tilde\rho g_{\eta},
\end{equation}
where $\tilde\rho=r^{2\eta+N-1}\exp\{\frac{r^2}{4D_v}\}$.
Integrating \eqref{eq4-16} from $r_1$ to $r$, and using \eqref{eq4-14} yields
\begin{align}
D_v\tilde\rho(r) g_{\eta}'(r)&\ge D_v\tilde\rho(r_1) g_{\eta}'(r_1)+\frac12\left(\kappa-\frac{\eta}2\right)\int_{r_1}^r\tilde\rho(s) g_{\eta}(s)ds
\nonumber
\\
\label{eq4-16}
&\ge D_v\tilde\rho(r_1) g_{\eta}'(r_1)+\frac{C^*}2\left(\kappa-\frac{\eta}2\right)\int_{r_1}^r s^{\eta-1}\exp\{\frac{s^2}{4D_v}\}ds.
\end{align}
From the above inequality, it is easy to obtain that there exists $r_2>r_1$ such that $\tilde\rho(r) g_{\eta}'(r)>0$ for any $r\ge r_2$ since
$$
\int_{r_1}^r s^{\eta-1}\exp\left\{\frac{s^2}{4D_v}\right\}ds\to+\infty, \  \text{as $r\to +\infty$}.
$$
Hence,
$$
g_{\eta}(r)>g_{\eta}(r_2)>0.
$$
That is
$$
\phi(r)>Cr^{\eta}
$$
for large $r$. We complete the proof of this lemma. \hfill $\Box$

\medskip

From Lemma \ref{lem4-1} and Lemma \ref{addlem4-3}, it is evident that for any small positive constant $\varepsilon$,
$$
C_1r^{2\kappa-\varepsilon}<\phi(r)<C_2r^{2\kappa+\varepsilon}
$$
when $r$ is sufficiently large. Additionally, Lemma \ref{lem4-2} indicates that when $\kappa\ge 0$,
$$
\phi(r)\ge C_1r^{2\kappa}
$$
for large $r$.
Next, we will show that when $-\frac N2\le \kappa<0$, the inequality $\phi(r)\ge C_1r^{2\kappa}$ also holds for  large $r$.
\begin{lemma}\label{lem4-4}
Assume that $-\frac N2\le \kappa< 0$, and $\alpha+\frac{\chi}{D_u}\ge 1$.
Let $\phi(r)\in C^2[0,\infty)$ be the global positive solution of \eqref{1-8}. Then when $AB^{\alpha-1}$ is appropriately small, there exists a positive constant $C^*$ such that
$$
\phi(r)\ge C^* r^{2\kappa},
$$
for large $r$.
In particular, when $\alpha+\frac{\chi}{D_u}=1$, the above estimate holds for any $A, B>0$.
\end{lemma}

{\it\bfseries Proof.} Recall Lemma \ref{lem4-1},  for any $0>\eta>2\kappa$, we have $g_{\eta}'(r)< 0$ when $r$ is large.
Consequently,
$$
\phi'(r)=(r^{\eta}g_\eta(r))'=\eta r^{\eta-1}g_\eta(r)+r^{\eta}g_\eta'(r)<0
$$
for large $r$.  By invoking Lemma \ref{lem4-1} and Lemma \ref{addlem4-3}, we find that for any small positive constant $\varepsilon$,
$$
C_1r^{2\kappa-\varepsilon}<\phi(r)<C_2r^{2\kappa+\varepsilon}
$$
when $r$ is large. In summary,  there exists $r_2>1$ such that
\begin{equation}\label{eq4-17}
C_1r^{2\kappa-\frac18}<\phi(r)<C_2r^{2\kappa+\frac 18}, \ \text{and}\ \phi'(r)<0, \text {when $r\ge r_2$},
\end{equation}
where $C_1$, $C_2$ are positive constants.
Next, direct integration of \eqref{1-9} from $r_2$ to $r$ gives
$$
D_v \rho(r)\phi'(r)=D_v \rho(r_2)\phi'(r_2)+\int_{r_2}^r\left(\kappa \rho\phi+A B^{\alpha-1}\rho e^{-\frac{r^2}{4D_u}} \phi^{\alpha+\frac{\chi}{D_u}}\right) ds>\int_{r_2}^r \kappa \rho(s)\phi(s) ds+D_v \rho(r_2)\phi'(r_2).
$$
Observing that $\kappa<0$, we subsequently apply \eqref{eq4-17} to further derive
\begin{align}
\frac{\phi'(r)}{\phi(r)}>\frac{\kappa}{D_v}\frac{1}{\rho(r)\phi(r)}\int_{r_2}^r  \rho(s)\phi(s) ds+\frac{\rho(r_2)\phi'(r_2)}{\rho(r)\phi(r)}\nonumber
\\
\label{eq4-18}
> \frac{C_2\kappa}{C_1D_v}\frac{1}{\rho(r)r^{2\kappa-\frac18}}\int_{r_2}^r\rho(s)s^{2\kappa+\frac 18} ds+\frac{\rho(r_2)\phi'(r_2)}{C_1r^{2\kappa-\frac18}\rho(r)}
\end{align}
when $r>r_2$. It is easy to verify that
$$
\lim_{r\to\infty}\frac{1}{r^{-\frac34}\rho(r)r^{2\kappa-\frac18}}\int_{r_2}^r\rho(s)s^{2\kappa+\frac 18} ds=2D_v.
$$
Hence, there exists $r_3>r_2$ such that when $r>r_3$,
$$
D_v<\frac{1}{r^{-\frac34}\rho(r)r^{2\kappa-\frac18}}\int_{r_2}^r\rho(s)s^{2\kappa+\frac 18} ds<3D_v.
$$
Substituting the above inequality  into \eqref{eq4-18},  yields
\begin{align*}
0>\frac{\phi'(r)}{\phi(r)}
> \frac{3C_2\kappa}{C_1}r^{-\frac34}+\frac{\rho(r_2)\phi'(r_2)}{C_1r^{2\kappa-\frac18}\rho(r)}>
-A^* r^{-\frac34}, \ \text{when $r>r_3$},
\end{align*}
where $A^*$ is a positive constant.

The proof that follows bears resemblance to the proof presented for Lemma \ref{lem4-2}.
To proceed, we define a function
we introduce a function defined as in Lemma \ref{lem4-2}, that is, we let $h(r)=\ln{\phi(r)}$, which satisfies the equation \eqref{eq4-5}.
Consequently,
\begin{align*}
D_v(\rho(r)h')'=-D_v|h'|^2\rho(r)+\kappa\rho(r)+AB^{\alpha-1}\rho(r)
e^{\left(\alpha+\frac{\chi}{D_u}-1\right)h}e^{-\frac{1}{4D_u}r^2}
>\kappa\rho(r)-\frac{D_vA^{*2}}{r^{\frac32}}\rho(r), \ \text{when $r>r_3$}.
\end{align*}
Integrating the above inequality from $r_3$ to $r$ gives
 \begin{align}
h(r)&>h(r_3)+\int_{r_3}^r\frac{\rho(r_3)h'(r_3)}{\rho(s)}ds+\frac{\kappa}{D_v}\int_{r_3}^r\frac{1}{\rho(\tau)}\int_{r_3}^{\tau}\rho(s)dsd\tau
-D_vA^{*2}\int_{r_3}^r\frac{1}{\rho(\tau)}\int_{r_3}^{\tau}s^{-\frac32}\rho(s)dsd\tau\nonumber
\\
\nonumber
&>\frac{\kappa}{D_v}\int_{r_3}^r\frac{1}{\rho(\tau)}\int_{r_3}^{\tau}\rho(s)dsd\tau-M_1
\\
\label{eq4-19}
&>2\kappa\ln r-M_2,
\end{align}
and this lemma is proved.
 \hfill $\Box$

\medskip

From the above two lemmas, it is not difficult to see that $\phi$ exhibits an algebraic growth or decay rate at infinity, and the growth rate (or decay rate for $\kappa<0$) is likely to be $r^{2\kappa}$. Below, we will prove this conclusion.
Take $\eta=2\kappa$ (namely, $\phi(r)=r^{2\kappa}g(r)$) in \eqref{4-1} to obtain
\begin{equation}
\label{4-11}
D_vg''+\left(\frac{r}2+\frac{D_v(4\kappa+N-1)}r\right)g'=-\frac{2\kappa D_v(2\kappa+N-2)}{r^2}g+AB^{\alpha-1}r^{2\kappa(\alpha+\frac{\chi}{D_u}-1)}g^{\alpha+\frac{\chi}{D_u}}e^{-\frac{r^2}{4D_u}}.
\end{equation}
We observe that the first term on the right-hand side of Equation  \eqref{4-11} will be the primary factor
since $g$ has an algebraic growth or decay rate. Therefore, the symbol of $\kappa(2\kappa+N-2)$  is essential.
To begin, we analyze the case $\kappa(2\kappa+N-2)>0$.

\begin{lemma}\label{lem4-5}
Assume that $\alpha+\frac{\chi}{D_u}\ge 1$, $\kappa\ge -\frac{N}2$, and $\kappa(2\kappa+N-2)>0$. Let $\phi(r)\in C^2[0,\infty)$ be the global positive solution of \eqref{1-8}.
Then there exists $M^*>0$ such that
$$
\lim_{r\to\infty} \frac{\phi(r)}{r^{2\kappa}}=M^*.
$$
\end{lemma}

{\it\bfseries Proof.} Recall the equation \eqref{4-11}, and observe that, under the assumptions of this lemma,
$2\kappa D_v(2\kappa+N-2)>0$. Furthermore,  from Lemma \ref{lem4-1}, it is evident that
$$
AB^{\alpha-1}r^{2\kappa(\alpha+\frac{\chi}{D_u}-1)}g^{\alpha+\frac{\chi}{D_u}}e^{-\frac{r^2}{4D_u}}<
\frac{\kappa D_v(2\kappa+N-2)}{r^2}g
$$
when $r$ is sufficiently large since $g<Cr^{\varepsilon}$ for some $\varepsilon>0$. It implies that there exists $R_0>0$ such that
\begin{align}
\label{4-12}
D_vg''+\left(\frac{r}2+\frac{D_v(4\kappa+N-1)}r\right)g'<-\frac{\kappa D_v(2\kappa+N-2)}{r^2}g,\quad \text{for $r>R_0$},
\end{align}
namely,
\begin{align}
\label{4-13}
D_v\left(r^{4\kappa+N-1}e^{\frac{r^2}{4D_v}}g'\right)'<-\kappa D_v(2\kappa+N-2)r^{4\kappa+N-3}e^{\frac{r^2}{4D_v}}g(r)<0,
\quad \text{for $r>R_0$}.
\end{align}
Therefore,
\begin{align}
\label{4-14}
r^{4\kappa+N-1}e^{\frac{r^2}{4D_v}}g'(r)<R_0^{4\kappa+N-1}e^{\frac{R_0^2}{4D_v}}g'(R_0).
\end{align}
If $g'(r)>0$ for any $r\ge R_0$,
we further have
$$
g(r)<g(R_0)+R_0^{4\kappa+N-1}e^{\frac{R_0^2}{4D_v}}g'(R_0)\int_{R_0}^rs^{-4\kappa-N+1}e^{-\frac{s^2}{4D_v}}ds,
$$
which means that $g$ is bounded.

In what follows, we show that $g'(r)>0$ for $r\ge R_0$  is impossible. In fact, if it does occur, it can be seen from the above inequality that $g$ has an upper bound, so there exists $M^*>0$ that such that
\begin{equation}\label{4-15}
g(r)\nearrow M^*, \ \text{as $r\to\infty$}.
\end{equation}
From \eqref{4-14}, it is easy to see that
$$
r^{L}g'(r)\to 0, \text{ as $r\to\infty$ for any $L\ge 0$}.
$$
Recalling \eqref{4-11}, we see that
$$
D_vr^2g''+\left(\frac{r^3}2+D_v(4\kappa+N-1)r\right)g'=-2\kappa D_v(2\kappa+N-2)g+AB^{\alpha-1}r^{2+2\kappa(\alpha+\frac{\chi}{D_u}-1)}g^{\alpha+\frac{\chi}{D_u}}e^{-\frac{r^2}{4D_u}}.
$$
Letting $r\to\infty$ gives
$$
\lim_{r\to\infty}r^2g''(r)=-2\kappa (2\kappa+N-2)M^*<0.
$$
While, on the other hands, by L'Hospital principle, we also have
$$
0=\lim_{r\to\infty} rg'(r)=\lim_{r\to\infty} \frac{g'(r)}{\frac{1}r}=\lim_{r\to\infty} \frac{g''(r)}{-\frac{1}{r^{2}}}=-\lim_{r\to\infty}r^{2}g''(r).
$$
It is a contradiction.

Thus, there exists $R_1\ge R_0$ such that $g'(R_1)\le 0$, from \eqref{4-13}, we have $g'(r)<0$ for any $r>R_1$.
Therefore,
$\displaystyle
\lim_{r\to\infty} g(r)
$
exists, which together with Lemma \ref{lem4-2} to complete the proof. \hfill $\Box$

Next, we consider the case $\kappa (2\kappa+N-2)\le 0$, which together with Lemma \ref{lem4-5} gives

\begin{lemma}\label{lem4-6}
Assume that $\alpha+\frac{\chi}{D_u}\ge 1$, $\kappa\ge -\frac{N}2$, and $\kappa (2\kappa+N-2)\le 0$. Let $\phi(r)\in C^2[0,\infty)$ be the global positive solution of \eqref{1-8}. Then
$$
\lim_{r\to\infty}\frac{\phi(r)}{r^{2\kappa}}=M^*,
$$
where $M^*$ is a positive constant.
\end{lemma}

{\it\bfseries Proof.}
Let
$$
\ln g=f.
$$
Then
\begin{equation}
\label{4-16}
D_v\left(f''+\left|\frac{f'}{f}\right|^2\right)+\left(\frac{r}2+\frac{D_v(4\kappa+N-1)}r\right)f'=\frac{2\kappa D_v(2-2\kappa-N)}{r^2}+AB^{\alpha-1}r^{2\kappa(\alpha+\frac{\chi}{D_u}-1)}g^{\alpha+\frac{\chi}{D_u}-1}e^{-\frac{r^2}{4D_u}}.
\end{equation}
which is equivalent to
\begin{align}
&D_v\left(f'r^{4\kappa+N-1}e^{\frac{r^2}{4D_v}}\right)'+
D_vr^{4\kappa+N-1}e^{\frac{r^2}{4D_v}}\left|\frac{f'}{f}\right|^2 \nonumber
\\
\label{4-17}
&=2\kappa D_v(2-2\kappa-N)r^{4\kappa+N-3}e^{\frac{r^2}{4D_v}}
+AB^{\alpha-1}e^{\frac{r^2}{4D_v}}r^{4\kappa+N-1+2\kappa(\alpha+\frac{\chi}{D_u}-1)}g^{\alpha+\frac{\chi}{D_u}-1}e^{-\frac{r^2}{4D_u}}.
\end{align}
Noticing that there exists $r_0>0$ such that
$$
\frac{2\kappa D_v(2-2\kappa-N)}{r^2}+AB^{\alpha-1}r^{2\kappa(\alpha+\frac{\chi}{D_u}-1)}g^{\alpha+\frac{\chi}{D_u}-1}e^{-\frac{r^2}{4D_u}}\le
\frac{CD_v}{r^2}
$$
for large $r>r_0$ since $g$ has an algebraic growth or decay rate, where $C$ is a positive constant. Substituting the above inequality into \eqref{4-17} yields
\begin{align}\label{4-18}
\left(f'r^{4\kappa+N-1}e^{\frac{r^2}{4D_v}}\right)'\le C r^{4\kappa+N-3}e^{\frac{r^2}{4D_v}}.
\end{align}
Integrating \eqref{4-18} from $r_0$ to $r$ twice to obtain
\begin{align}\label{4-19}
f(r)\le f(r_0)+f'(r_0)r_0^{4\kappa+N-1}e^{\frac{r_0^2}{4D_v}}\int_{r_0}^rs^{-4\kappa-N+1}e^{-\frac{s^2}{4D_v}}ds
+ C \int_{r_0}^rs^{-4\kappa-N+1}e^{-\frac{s^2}{4D_v}}  \int_{r_0}^s\tau^{4\kappa+N-3}e^{\frac{\tau^2}{4D_v}}d\tau ds.
\end{align}
Similar to the proof of Lemma \ref{lem2-3} or Lemma \ref{lem4-2} to get that
\begin{align*}
s^{-4\kappa-N+1}e^{-\frac{s^2}{4D_v}}  \int_{r_0}^s\tau^{4\kappa+N-3}e^{\frac{\tau^2}{4D_v}}d\tau\le 2D_v\frac{1}{s^3}
+4D_v^2|4\kappa+N-4|\frac{1}{s^4},
\end{align*}
it implies that $f(r)$ has an upper bound, i.e. $g$ has an upper bound. That is there exists $C^{**}>0$ such that
\begin{equation}\label{4-20}
\frac{\phi(r)}{r^{2\kappa}}< C^{**}
\end{equation}
for $r>r_0$.
From \eqref{4-11}, we also see that $g(r)$ is  monotonic for large $r$ since the right side of the equal sign is always greater than $0$.
In fact, if there exists $r_1>0$ such that $g'(r_1)\ge 0$, then from \eqref{4-11} it is easy to get that $g'(r)>0$ for $r>r_1$.
Then either $g$ strictly monotonically decreases, or there is $R$ such that $g$ monotonically increase when $r>R$, which together with Lemma
\ref{lem4-2} and \eqref{4-20}, we obtain that there exists $\tilde C>0$ such that
$$
\lim_{r\to\infty}\frac{\phi(r)}{r^{2\kappa}}=\tilde C.
$$
We complete the proof of this lemma.
\hfill $\Box$

Combining Lemma \ref{lem4-5} and Lemma \ref{lem4-6}, we complete the proof of Theorem \ref{thm2-2}.

 In summary, we have accomplished the asymptotic rate estimation of $\phi$ as $r\to\infty$, thereby concluding the proof of Theorem \ref{thm2-2}.

 Recalling \eqref{1-6}, we see that
$$
\varphi(r)=A \phi^{\frac{\chi}{D_u}}e^{-\frac{1}{4D_u}r^2}.
$$
Therefore,
$$
u(x,t)=A t^{-\frac{N}2}|\phi(r)|^{\frac{\chi}{D_u}}e^{-\frac{1}{4D_u}r^2}, \qquad v(x,t)=Bt^{\kappa}\phi(r),
$$
with $r=t^{-\frac12}|x|$. A direct calculation yields
$$
\int_{\mathbb R^N}u(x,t)dx=\int_{\mathbb R^N}A |\phi(t^{-\frac12}|x|)|^{\frac{\chi}{D_u}}e^{-\frac{|x|^2}{4D_u t}}dt^{-\frac12}x
=
|\partial B_1|\int_0^\infty A  |\phi(r)|^{\frac{\chi}{D_u}}e^{-\frac{1}{4D_u}r^2}r^{N-1}dr.
$$
Notice that $\phi\in C^2[0,\infty)$, and as $r\to\infty$, $\phi\sim r^{2\kappa}$ for $\kappa\ge -\frac{N}2$; $|\phi|=o(\rho^{-\frac12})$ for $\kappa< -\frac{N}2$. Hence $\int_0^\infty A  |\phi(r)|^{\frac{\chi}{D_u}}e^{-\frac{1}{4D_u}r^2}r^{N-1}dr$ is integrable. Let
\begin{equation}\label{4-21}
M:=\int_{\mathbb R^N}u(x,t)dx=|\partial B_1|\int_0^\infty A  |\phi(r)|^{\frac{\chi}{D_u}}e^{-\frac{1}{4D_u}r^2}r^{N-1}dr.
\end{equation}
Next, we give following proposition.

\begin{proposition}\label{pro-2}
Assume that $(1-\alpha)\kappa=1-\frac N2$. Let $(u, v)$ be the self-similar solution of the model \eqref{1-1} obtained in the above three theorems, which is defined as in \eqref{1-2}. That is,
$$
u(x,t)=A t^{-\frac{N}2}|\phi(r)|^{\frac{\chi}{D_u}}e^{-\frac{1}{4D_u}r^2}, \qquad v(x,t)=Bt^{\kappa}\phi(r),
$$
with $r=t^{-\frac12}|x|$. It is easy to see that both $u$ and $v$ are smooth when $t>0$, and $u$ is bounded for any $t>0$.

Denote
$$
M:=\int_{\mathbb R^N}u(x,t)dx.
$$
Then for any $\kappa\in \mathbb R$,
$$
u(x,t)\to M\delta(x), \quad \text{as}\ t\to 0^+
$$
in the sense of distribution. And for any $1<p\le \infty$
\begin{equation}\label{1-11}
\lim_{t\to\infty}t^{\frac{N}2(1-\frac1p)}\|u(\cdot, t)\|_{L^p}=M_p,
\end{equation}
where $M_p$ are positive constants depending only on $p$.

Regarding the singularity of $v$ at the initial moment, it varies depending on the value of $\kappa$, which can be regular, very singular, or less singular. We present the results in the following three cases.

i) when $\kappa\ge 0$, $v$ is regular.

ii) When $\kappa\le -\frac{N}{2}$, $v$ is very singular. More precisely, when $\kappa<-\frac{N}{2}$ ($\Leftrightarrow -2\kappa-N>0$),
$$
|x|^{-2\kappa-N}v(x,t)\to M_1\delta(x), \quad \text{as}\ t\to 0^+,
$$
in the sense of distribution, where
$$
M_1=B|\partial B_1|\int_0^\infty r^{-2\kappa-1}\phi(r)dr;
$$
when $\kappa=-\frac N2$, for any given $\varepsilon>0$,
$$
\lim_{t\to 0^+}\int_{|x|<\varepsilon}v(x,t)dx=\infty.
$$

iii) When $-\frac{N}{2}< \kappa<0$,  $v$ is less singular. That is
$$
\lim_{t\to 0^+}\int_{|x|<\varepsilon}v(x,t)dx\le C\varepsilon^{2\kappa+N}, \ \text{for $-\frac N2<\kappa<0$}.
$$
While, for any $p>\frac{N}{2|\kappa|}$,
$$
\lim_{t\to 0^+}\int_{|x|<\varepsilon}|v(x,t)|^pdx\to \infty, \ \text{for $-\frac N2<\kappa<0$ with $p>\frac{N}{2|\kappa|}$.}
$$

iv) When $\kappa<0$,  $v$ algebraically decays to $0$ in the sense of $L^p$-norm for some $p\ge 1$ as $t\to\infty$.
More precisely, when $\kappa<-\frac N2$, for any $1\le p\le \infty$, we have  $-\kappa-\frac N{2p}>\frac N2(1-\frac 1p)$,
and
$$
\lim_{t\to\infty} t^{-\kappa-\frac N{2p}}\|v(\cdot, t)\|_{L^p}=\tilde M_p,
$$
with $\tilde M_p=\left(\int_0^{\infty}|\phi(r)|^p r^{N-1}dx\right)^{\frac 1p}$ is a positive constant depending on $p$;
when $-\frac N2\le \kappa<0$, take $p^*=\frac{N}{2|\kappa|}$, then when  $p>p^*$, we have
$$
\lim_{t\to\infty} t^{-\kappa-\frac N{2p}}\|v(\cdot, t)\|_{L^p}=\tilde M_p,
$$
with  $-\kappa-\frac N{2p}>0$.
\end{proposition}

{\it\bfseries Proof.} Note that for any $f\in C_0(\mathbb R^N)$,
\begin{align}\label{4-22}
\int_{\mathbb R^N}u(x,t)f(x)dx=\int_{\mathbb R^N}A t^{-\frac{N}2}|\phi(t^{-\frac12}|x|)|^{\frac{\chi}{D_u}}e^{-\frac{|x|^2}{4D_u t}}
f(x)dx=\int_{\mathbb R^N}A|\phi(|\xi|)|^{\frac{\chi}{D_u}}e^{-\frac{|\xi|^2}{4D_u}}f(\sqrt t \xi)d\xi.
\end{align}
Then
\begin{align}
&\left|\int_{\mathbb R^N}u(x,t)f(x)dx-Mf(0)\right|=\left|\int_{\mathbb R^N}A|\phi(|\xi|)|^{\frac{\chi}{D_u}}e^{-\frac{|\xi|^2}{4D_u}}\Big(f(\sqrt t \xi)-f(0)\Big)d\xi\right|\nonumber
\\
\le &\left|\int_{|\xi|\le R}A|\phi(|\xi|)|^{\frac{\chi}{D_u}}e^{-\frac{|\xi|^2}{4D_u}}\Big|f(\sqrt t \xi)-f(0)\Big|d\xi\right|
+\left|\int_{|\xi|> R}A|\phi(|\xi|)|^{\frac{\chi}{D_u}}e^{-\frac{|\xi|^2}{4D_u}}\Big|f(\sqrt t \xi)-f(0)\Big|d\xi\right|\nonumber
\\
\label{4-23}
\le &M\sup_{|\xi|\le R}\Big|f(\sqrt t \xi)-f(0)\Big|
+2\|f\|_{L^\infty}\left|\int_{|\xi|> R}A|\phi(|\xi|)|^{\frac{\chi}{D_u}}e^{-\frac{|\xi|^2}{4D_u}}d\xi\right|.
\end{align}
Recalling \eqref{4-21}, for any $\varepsilon>0$, there exists $R_0>0$ such that
\begin{equation}\label{4-24}
|\partial B_1|\int_{R_0}^\infty A  |\phi(r)|^{\frac{\chi}{D_u}}e^{-\frac{1}{4D_u}r^2}r^{N-1}dr<\frac{\varepsilon}{4\|f\|_{L^\infty}}.
\end{equation}
For the above established $R_0$, there exists $t_0>0$ such that when $t<t_0$,
\begin{equation}\label{4-25}
\sup_{|\xi|\le R_0}\Big|f(\sqrt t \xi)-f(0)\Big| <\frac{\varepsilon}{2M}
\end{equation}
Combining \eqref{4-23}, \eqref{4-24} and \eqref{4-25}, we obtain that for any $\varepsilon>0$, there exists $t_0>0$, such that when $t<t_0$,
$$
\left|\int_{\mathbb R^N}u(x,t)f(x)dx-Mf(0)\right|<\varepsilon.
$$
It implies that
$$
u(x,t)\to M\delta(x), \quad \text{as}\ t\to 0^+.
$$
On the other hand, we note that for any $p\in(1,\infty)$,
\begin{align*}
t^{\frac{N(p-1)}2}\int_{\mathbb R^N}|u|^pdx&=t^{\frac{N(p-1)}2}\int_{\mathbb R^N}A^p t^{-\frac{Np}2}|\phi(t^{-\frac12}|x|)|^{\frac{p\chi}{D_u}}e^{-\frac{p|x|^2}{4D_u t}}dx
\\
&=|\partial B_1|\int_0^\infty A^p |\phi(r)|^{\frac{p\chi}{D_u}}e^{-\frac{pr^2}{4D_u}}r^{N-1}dr
\\
&=M_p^p,
\end{align*}
that is
\begin{equation}\label{Lp}
t^{\frac{N}2(1-\frac1p)}\|u(\cdot, t)\|_{L^p}\equiv M_p.
\end{equation}
when $p=\infty$, we also have
\begin{equation}\label{Linfty}
t^{\frac{N}2}\|u(\cdot, t)\|_{L^\infty}\equiv M_\infty.
\end{equation}

Next,  we examine the singularity of $v$ at the initial time. We discuss the following three cases separately: i) when $\kappa<-\frac{N}{2}$, ii) when $-\frac{N}{2}\le \kappa<0$, and iii) when $\kappa\ge 0$.

i) Firstly, when $\kappa<-\frac N2$, we see that
$$
\int_{\mathbb R^N}|x|^{-2\kappa-N}v(x,t)dx=B\int_{\mathbb R^N}|x|^{-2\kappa-N}t^{\kappa}\phi(t^{-\frac12}|x|)dx=B|\partial B_1|\int_0^\infty r^{-2\kappa-1}\phi(r)dr.
$$
It is clear that $\int_0^\infty r^{-2\kappa-1}\phi(r)dr<\infty$ since $\phi(r)=o(\rho^{-\frac{1}2})$.
Let
$$
M_1=B|\partial B_1|\int_0^\infty r^{-2\kappa-1}\phi(r)dr.
$$
Then for any $f\in C_0(\mathbb R^N)$,
\begin{align}\label{4-26}
\int_{\mathbb R^N}|x|^{-2\kappa-N}v(x,t)f(x)dx=\int_{\mathbb R^N}|x|^{-2\kappa-N}t^{\kappa}\phi(t^{-\frac12}|x|)f(x)dx
=\int_{\mathbb R^N}|\xi|^{-2\kappa-N}\phi(\xi)f(\sqrt t\xi)d\xi.
\end{align}
Hence,
\begin{align}
&\left|\int_{\mathbb R^N}|x|^{-2\kappa-N}v(x,t)f(x)dx-Mf(0)\right|=\left|\int_{\mathbb R^N}|\xi|^{-2\kappa-N}\phi(\xi) \Big(f(\sqrt t \xi)-f(0)\Big)d\xi\right|\nonumber
\\
\le &\left|\int_{|\xi|\le R}|\xi|^{-2\kappa-N}\phi(\xi)\Big|f(\sqrt t \xi)-f(0)\Big|d\xi\right|
+\left|\int_{|\xi|> R}|\xi|^{-2\kappa-N}\phi(\xi)\Big|f(\sqrt t \xi)-f(0)\Big|d\xi\right|\nonumber
\\
\label{4-27}
\le &M_1\sup_{|\xi|\le R}\Big|f(\sqrt t \xi)-f(0)\Big|
+2\|f\|_{L^\infty}\left|\int_{|\xi|> R}|\xi|^{-2\kappa-N}\phi(\xi)d\xi\right|.
\end{align}
Similar to \eqref{4-24} and \eqref{4-25}, we finally obtain that for any $\varepsilon>0$, there exists $t_0>0$, such that when $t<t_0$,
$$
\left|\int_{\mathbb R^N}|x|^{-2\kappa-N}v(x,t)f(x)dx-M_1f(0)\right|<\varepsilon.
$$
It implies that
$$
|x|^{-2\kappa-N}v(x,t)\to M_1\delta(x), \quad \text{as}\ t\to 0^+.
$$

ii) When $-\frac{N}{2}\le \kappa<0$. By Theorem \ref{thm2-2}, there exists $R>0$ such that
$$
\frac{M^*}2r^{2\kappa}\le \phi(r)\le \frac{3M^*}2r^{2\kappa}, \ \text{ for $r\ge R$},
$$
it implies that
$$
\frac{M^*}2|x|^{2\kappa}\le v(x,t)=t^{\kappa}\phi\left(\frac{|x|}{\sqrt t}\right)\le  \frac{3M^*}2|x|^{2\kappa}, \ \text{ for $|x|\ge R\sqrt t$}.
$$
When $\kappa=-\frac N2$, for any given $\varepsilon>0$,
\begin{align*}
\int_{|x|<\varepsilon}v(x,t)dx=&t^{-\frac N2}\int_{|x|<R\sqrt t}\phi\left(\frac{|x|}{\sqrt t}\right)dx+
t^{-\frac N2}\int_{R\sqrt t\le|x|<\varepsilon}\phi\left(\frac{|x|}{\sqrt t}\right)dx
\\
\ge & |\partial B_1|\int_0^R\phi(r)r^{N-1}dr+\frac{M^*}2
\int_{R\sqrt t\le|x|<\varepsilon}|x|^{-N}dx
\\
\ge &\frac{M^*}2|\partial B_1|
\int_{R\sqrt t}^{\varepsilon}r^{-1}dr
\\
\ge &\frac{M^*}2|\partial B_1|\ln\left(\frac{\varepsilon}{R\sqrt t}\right).
\end{align*}
Hence, for any given $\varepsilon>0$,
$$
\lim_{t\to 0^+}\int_{|x|<\varepsilon}v(x,t)dx=\infty,  \ \text{for $\kappa=-\frac N2$}.
$$
When $-\frac N2<\kappa<0$,
\begin{align*}
\int_{|x|<\varepsilon}v(x,t)dx=&t^{\kappa}\int_{|x|<R\sqrt t}\phi\left(\frac{|x|}{\sqrt t}\right)dx+
t^{\kappa}\int_{R\sqrt t\le|x|<\varepsilon}\phi\left(\frac{|x|}{\sqrt t}\right)dx
\\
\le &t^{\kappa+\frac N2}|\partial B_1|\int_0^R\phi(r)r^{N-1}dr+\frac{3M^*}2
\int_{R\sqrt t\le|x|<\varepsilon}|x|^{2\kappa}dx
\\
\le &Ct^{\kappa+\frac N2}+C\varepsilon^{2\kappa+N},
\end{align*}
thus,
\begin{equation}\label{4-28}
\lim_{t\to 0^+}\int_{|x|<\varepsilon}v(x,t)dx\le C\varepsilon^{2\kappa+N}, \ \text{for $-\frac N2<\kappa<0$}.
\end{equation}
While, by taking $p>\frac{N}{2|\kappa|}$ yields
\begin{align*}
\int_{|x|<\varepsilon}|v(x,t)|^pdx=&t^{\kappa p}\int_{|x|<R\sqrt t}\phi^p\left(\frac{|x|}{\sqrt t}\right)dx+
t^{\kappa p}\int_{R\sqrt t\le|x|<\varepsilon}\phi^p\left(\frac{|x|}{\sqrt t}\right)dx
\\
\ge &\left(\frac{M^*}2\right)^p|\partial B_1|
\int_{R\sqrt t}^{\varepsilon}r^{2\kappa p+N-1}dr \to \infty \text{as $t\to 0^+$},
\end{align*}
that is,
\begin{equation}\label{4-29}
\lim_{t\to 0^+}\int_{|x|<\varepsilon}|v(x,t)|^pdx\to \infty, \ \text{for $-\frac N2<\kappa<0$ with $p>\frac{N}{2|\kappa|}$.}
\end{equation}

iii) While if $\kappa\ge 0$, from Theorem \ref{thm2-2}, it is easy to see that $v$  is non-singular as $t\to 0^+$.

iv) At last, we consider the large time behavior. We also note that
\begin{align*}
t^{-\kappa p-\frac N2}\int_{\mathbb R^N}|v|^pdx=t^{-\kappa p-\frac N2}\int_{\mathbb R^N}t^{\kappa p}\left|\phi(\frac{|x|}{\sqrt t})\right|^pdx
=\int_0^{\infty}|\phi(r)|^p r^{N-1}dx.
\end{align*}
It is easy to see that when $\kappa<-\frac N2$, we have  $-\kappa p-\frac N2>\frac N2(p-1)$, and $\int_0^{\infty}|\phi(r)|^p r^{N-1}dx$ is integrable for any $p\ge 1$;
when $0>\kappa\ge -\frac N2$, $\int_0^{\infty}|\phi(r)|^p r^{N-1}dx$ is integrable when $2\kappa p+N<0$.

We complete the proof.
\hfill $\Box$

Theorem \ref{thm2-4} is a direct result of Proposition \ref{pro-2}.
\medskip

\medskip
{\bf Declarations:} The author declares that there is no competing interest.

\medskip
{\bf Data availability statement:} The author confirms that the data supporting the findings of this
study are available within the article

\end{document}